\documentclass[twoside,12pt]{article}
\usepackage{amssymb}

\usepackage{a4wide}
\usepackage{amsthm}
\usepackage{amsmath}
\usepackage[all]{xy}
\usepackage{enumerate}
\usepackage{fancyhdr}

\newtheorem{theorem}{Theorem}[section]
\newtheorem{proposition}[theorem]{Proposition}
\newtheorem{corollary}[theorem]{Corollary}
\newtheorem{lemma}[theorem]{Lemma}
\theoremstyle{definition}

\newtheorem*{Beweis}{Proof}
\newtheorem{definition}[theorem]{Definition}
\newtheorem{punto}[theorem]{}
\theoremstyle{remark}
\newtheorem{remark}[theorem]{Remark}
\newtheorem{remarks}[theorem]{Remarks}
\CompileMatrices
\swapnumbers
\CompileMatrices
\input{tcilatex}

\begin{document}

\title{Morita Contexts for Corings and Equivalences\thanks{%
MSC (2000): 16D90, 16S40, 16W30, 13B02 \newline
Keywords: Morita Contexts, Hopf Algebras, (Cleft) Galois Extensions,
Entwining Structures, Entwined Modules.}}
\author{\textbf{Jawad Y. Abuhlail}\thanks{%
Current Address: Department of Mathematical Sciences, Box \# 5046, KFUPM,
31261 Dhahran (Saudi Arabia), \textbf{email:} abuhlail@kfupm.edu.sa} \\
Mathematics Department, Birzeit University\\
Birzeit - Palestine}
\date{}
\maketitle

\begin{abstract}
In this note we study Morita contexts and Galois extensions for corings. For
a coring $\mathcal{C}$ over a (not necessarily commutative) ground ring $A$
we give equivalent conditions for $\mathcal{M}^{\mathcal{C}}$ to satisfy the
weak. resp. the strong structure theorem. We also characterize the so called 
\emph{cleft }$C$\emph{-Galois extensions} over commutative rings. Our
approach is similar to that of Y. Doi and A. Masuoka in their work on
(cleft) $H$-Galois extensions (e.g. \cite{Doi94}, \cite{DM92}).
\end{abstract}

\section*{Introduction}

Let $\mathcal{C}$ be a coring over a not necessarily commutative ring $A$
and assume $A$ to be a right $\mathcal{C}$-comodule through $\varrho
_{A}:A\longrightarrow A\otimes _{A}\mathcal{C}\simeq \mathcal{C},$ $a\mapsto 
\mathbf{x}a$ for some group-like element $\mathbf{x}\in \mathcal{C}$ (see 
\cite[Lemma 5.1]{Brz02}). In the first section we study from the viewpoint
of Morita theory the relationship between $A$ and its subring of
coinvariants $B:=A^{co\mathcal{C}}:=\{b\in A\mid \varrho (b)=b\mathbf{x}\}.$
We consider the $A$-ring $^{\ast }\mathcal{C}:=\mathrm{Hom}_{A-}(\mathcal{C}%
,A)$ and its left ideal $Q:=\{q\in $ $^{\ast }\mathcal{C}\mid \sum
c_{1}q(c_{2})=q(c)\mathbf{x}$ for all $c\in \mathcal{C}\}$ and show that $B$
and $^{\ast }\mathcal{C}$ are connected via a Morita context using $%
_{B}A_{^{\ast }\mathcal{C}}$ and $_{^{\ast }\mathcal{C}}Q_{B}$ as connecting
bimodules. Our Morita context is in fact a generalization of Doi's Morita
context presented in \cite{Doi94}.

In the second section we introduce the weak (resp. the strong) structure
theorem for $\mathcal{M}^{\mathcal{C}}.$ For the case $_{A}\mathcal{C}$ is
locally projective, in the sense of B. Zimmermann-Huignes,\ we characterize $%
A$ being a generator (a progenerator) in the category of right $\mathcal{C}$%
-comodules by $\mathcal{M}^{\mathcal{C}}$ satisfying the weak (resp. the
strong) structure theorem. Here the notion of Galois corings introduced by
T. Brzezi\'{n}ski \cite{Brz02} plays an important role. The results and
proofs are essentially module theoretic and similar to those of \cite{MZ97}
for the catgeory $\mathcal{M}(H)_{A}^{C}$ of Doi-Koppinen modules
corresponding to a right-right Doi-Koppinen structure $(H,A,C)$ (see also 
\cite{MSTW01} for the case $C=H$).

The notion of a $C$-Galois extension $A$ of a ring $B$ was introduced by T.
Brzezi\'{n}ski and S. Majid in \cite{BM98} and is related to the so called
entwining structures introduced in the same paper. In the third section we
give equivalent conditions for a $C$-Galois extension $A/B$ to be cleft. Our
results generalize results of \cite{Brz99} from the case of a base field to
the case of a commutative ground ring. In the special case $\varrho (a)=\sum
a_{\psi }\otimes x^{\psi },$ for some group-like element $x\in C,$ we get a
complete generalization of \cite[Theorem 1.5 ]{DM92} (and \cite[Theorem 2.5]
{Doi94}).

With $A$ we denote a not necessarily commutative ring with $1_{A}\neq 0_{A}$
and with $\mathcal{M}_{A}$ (resp. $_{A}\mathcal{M},$ $_{A}\mathcal{M}_{A}$)
the category of \emph{unital }right $A$-modules (resp. left $A$-modules, $A$%
-bimodules). For every right $A$-module $W$ we denote by $\mathrm{Gen}%
(W_{A}) $ (resp. $\sigma \lbrack W_{A}]$) the class of $W$-generated (resp. $%
W$-subgenerated) right $A$-modules. For the well developed theory of
categories of type $\sigma \lbrack W]$ the reader is referred to 
\cite[Section 15]{Wis88}.

An $A$-module $W$ is called \textbf{locally projective }(in the sense of B.
Zimmermann-Huignes \cite{Z-H76}), if for every diagram 
\begin{equation*}
\xymatrix{0 \ar[r] & F \ar@{.>}[dr]_{g' \circ \iota} \ar[r]^{\iota} & W
\ar[dr]^{g} \ar@{.>}[d]^{g'} & & \\ & & L \ar[r]_{\pi} & N \ar[r] & 0}
\end{equation*}
with exact rows and $F$ f.g.: for every $A$-linear map $g:W\longrightarrow
N, $ there exists an $A$-linear map $g^{\prime }:W\longrightarrow L$, such
that the entstanding parallelogram is commutative. Note that every
projective $A$-module is locally projective. By \cite[Theorem 2.1]{Z-H76} a
left $A$-module $W$ is locally projective, iff for every right $A$-module $M$
the following map is injective 
\begin{equation*}
\alpha _{M}^{W}:M\otimes _{A}W\longrightarrow \mathrm{Hom}_{-A}(^{\ast }W,M),%
\text{ }m\otimes _{A}w\mapsto \lbrack f\mapsto mf(w)].
\end{equation*}
It's easy then to see that every locally projective $A$-module is flat and $%
A $-cogenerated.

Let $\mathcal{C}$ be an $A$-coring. We consider the canonical $A$-bimodule $%
^{\ast }\mathcal{C}:=\mathrm{Hom}_{A-}(\mathcal{C},A)$ as an $A$-ring with
the canonical $A$-bimodule structure, multiplication $(f\cdot g)(c):=\sum
g(c_{1}f(c_{2}))$ and unity $\varepsilon _{\mathcal{C}}.$ If $_{A}\mathcal{C}
$ is locally projective, then we have an isomorphism of categories $\mathcal{%
M}^{\mathcal{C}}\simeq \sigma \lbrack \mathcal{C}_{^{\ast }\mathcal{C}}]$
(in particular $\mathcal{M}^{\mathcal{C}}\subseteq \mathcal{M}_{^{\ast }%
\mathcal{C}}$ is a full subcategory) and we have a left exact functor $%
\mathrm{Rat}^{\mathcal{C}}(-):\mathcal{M}_{^{\ast }\mathcal{C}}\rightarrow 
\mathcal{M}^{\mathcal{C}}$ assigning to every right $^{\ast }\mathcal{C}$%
-module its maximum $\mathcal{C}$\emph{-rational} $^{\ast }\mathcal{C}$%
-submodule, which turns to be a right $\mathcal{C}$-comodule. Moreover $%
\mathcal{M}^{\mathcal{C}}=\mathcal{M}_{^{\ast }\mathcal{C}}$ iff $_{A}%
\mathcal{C}$ is f.g. and projective. For more investigation of the $\mathcal{%
C}$-rational $^{\ast }\mathcal{C}$-modules see \cite{Abu03}.

After this paper was finished, it turned out that some results in this paper
were discovered independently by S. Caenepeel, J. Vercruysse and S. Wang 
\cite{CVW}.

\section{Morita Contexts}

In this section we fix the following: $\mathcal{C}$ is an $A$-coring with
group-like element $\mathbf{x}\ $and $A$ is a right $\mathcal{C}$-comodule
with structure map 
\begin{equation*}
\varrho _{A}:A\longrightarrow A\otimes _{A}\mathcal{C}\simeq \mathcal{C},%
\text{ }a\mapsto \mathbf{x}a
\end{equation*}
(e.g. \cite[Lemma 5.1]{Brz02}), hence $A\in \mathcal{M}_{^{\ast }\mathcal{C}%
} $ with $a\leftharpoonup g:=\sum a_{<0>}g(a_{<1>})=g(\mathbf{x}a)$ for all $%
a\in A$ and $g\in $ $^{\ast }\mathcal{C}.$ For $M\in \mathcal{M}_{^{\ast }%
\mathcal{C}}$ put 
\begin{equation*}
M^{\mathbf{x}}:=\{m\in M\mid mg=mg(\mathbf{x})\text{ for all }g\in \text{ }%
^{\ast }\mathcal{C}\}.
\end{equation*}
In particular $A^{\mathbf{x}}:=\{a\in A\mid a\leftharpoonup g=ag(\mathbf{x})$
for all $g\in $ $^{\ast }\mathcal{C}\}\subset A$ is a subring. For $M\in 
\mathcal{M}^{\mathcal{C}}$ we set 
\begin{equation*}
M^{co\mathcal{C}}:=\{m\in M\mid \varrho (m)=m\otimes _{A}\mathbf{x}%
\}\subseteq M^{\mathbf{x}}.
\end{equation*}
Obviously $B:=A^{co\mathcal{C}}=\{b\in A\mid b\mathbf{x}=\mathbf{x}%
b\}\subseteq A^{\mathbf{x}}$ is a subring and $\varrho _{A}$ is $(B,A)$%
-bilinear. For $M\in \mathcal{M}^{\mathcal{C}}$ we have $M^{co\mathcal{C}%
}\in \mathcal{M}_{B}.$ Moreover we set 
\begin{equation*}
Q:=\{q\in \text{ }^{\ast }\mathcal{C}\mid \sum c_{1}q(c_{2})=q(c)\mathbf{x}%
\text{ for all }c\in \mathcal{C}\}\subseteq (^{\ast }\mathcal{C})^{\mathbf{x}%
}.
\end{equation*}

\begin{lemma}
\label{Psi-pr}

\begin{enumerate}
\item  For every right $^{\ast }\mathcal{C}$-module $M$ we have an
isomorphism of right $B$-modules 
\begin{equation*}
\omega _{M}:\mathrm{Hom}_{-^{\ast }\mathcal{C}}(A,M)\longrightarrow M^{%
\mathbf{x}},\text{ }f\mapsto f(1_{A})
\end{equation*}
with inverse $m\mapsto \lbrack a\mapsto ma].$

\item  Let $_{A}\mathcal{C}$ be locally projective. If $M\in \mathcal{M}^{%
\mathcal{C}},$ then $M^{co\mathcal{C}}=M^{\mathbf{x}}\simeq \mathrm{Hom}%
_{-^{\ast }\mathcal{C}}(A,M)=\mathrm{Hom}^{\mathcal{C}}(A,M).$ Hence 
\begin{equation*}
\Psi _{M}:M^{co\mathcal{C}}\otimes _{B}A\longrightarrow M,\text{ }m\otimes
_{B}a\mapsto ma
\end{equation*}
is surjective \emph{(}resp. injective, bijective\emph{)}, iff 
\begin{equation*}
\Psi _{M}^{\prime }:\mathrm{Hom}^{\mathcal{C}}(A,M)\otimes
_{B}A\longrightarrow M,\text{ }f\otimes _{B}a\mapsto f(a)
\end{equation*}
is surjective \emph{(}resp. injective, bijective\emph{)}.

\item  We have $\mathrm{Hom}_{-^{\ast }\mathcal{C}}(A,^{\ast }\mathcal{C}%
)\simeq (^{\ast }\mathcal{C})^{\mathbf{x}}.$ If moreover $_{A}\mathcal{C}$
is $A$-cogenerated \emph{(}resp. locally projective and $^{\Box }\mathcal{C}%
:=\mathrm{Rat}^{\mathcal{C}}(^{\ast }\mathcal{C}_{^{\ast }\mathcal{C}})$%
\emph{)}, then $Q=(^{\ast }\mathcal{C})^{\mathbf{x}}$ \emph{(}resp. $%
Q=(^{\Box }\mathcal{C})^{co\mathcal{C}}$\emph{)}.

\item  For every $M\in \mathcal{M}_{^{\ast }\mathcal{C}}$ \emph{(}resp. $%
M\in \mathcal{M}^{\mathcal{C}}$\emph{)} and all $m\in M,$ $q\in Q$ we have $%
mq\in M^{\mathbf{x}}$ \emph{(}resp. $mq\in M^{co\mathcal{C}}$\emph{)}.
\end{enumerate}
\end{lemma}

\begin{Beweis}
\begin{enumerate}
\item  Obvious.

\item  Trivial.

\item  Considering $^{\ast }\mathcal{C}$ as a right $^{\ast }\mathcal{C}$%
-module via right multiplication we get $\mathrm{Hom}_{-^{\ast }\mathcal{C}%
}(A,^{\ast }\mathcal{C})\simeq (^{\ast }\mathcal{C})^{\mathbf{x}}$ by (1).\
If $q\in $ $(^{\ast }\mathcal{C})^{\mathbf{x}},$ then we have for all $g\in $
$^{\ast }\mathcal{C}$ and $c\in \mathcal{C}:$%
\begin{equation*}
g(\sum c_{1}q(c_{2}))=\sum g(c_{1}q(c_{2}))=(q\cdot g)(c)=(qg(\mathbf{x}%
))(c)=q(c)g(\mathbf{x})=g(q(c)\mathbf{x}),
\end{equation*}
i.e. $\sum c_{1}q(c_{2})-q(c)\mathbf{x}\in \mathrm{\func{Re}}(\mathcal{C},A)$
$:=\bigcap \{\mathrm{Ke}(g)\mid g\in \mathrm{Hom}_{A-}(\mathcal{C},A)\}.$ If 
$_{A}\mathcal{C}$ is $A$-cogenerated, then $\mathrm{\func{Re}}(\mathcal{C}%
,A)=0,$ hence $Q=(^{\ast }\mathcal{C})^{\mathbf{x}}.$

Assume $_{A}\mathcal{C}$ to be locally projective. Then we have for all $%
q\in Q,g\in $ $^{\ast }\mathcal{C}$ and $c\in \mathcal{C}:$%
\begin{equation*}
(q\cdot g)(c)=\sum g(c_{1}q(c_{2}))=g(q(c)\mathbf{x})=q(c)g(\mathbf{x})=(qg(%
\mathbf{x}))(c),
\end{equation*}
hence $q\in $ $^{\Box }\mathcal{C},$ with $\varrho (q)=q\otimes _{A}\mathbf{x%
},$ i.e. $q\in $ $(^{\Box }\mathcal{C})^{co\mathcal{C}}.$ On the other hand,
if $q\in (^{\Box }\mathcal{C})^{co\mathcal{C}},$ then for all $g\in $ $%
^{\ast }\mathcal{C}$ we have $q\cdot g=qg(\mathbf{x}),$ i.e. $q\in (^{\ast }%
\mathcal{C})^{\mathbf{x}}=Q.$

\item  Let $M\in \mathcal{M}_{^{\ast }\mathcal{C}}.$ Then we have for all $%
q\in Q,$ $g\in $ $^{\ast }\mathcal{C}$ and $m\in M:$%
\begin{equation*}
(mq)g=m(q\cdot g)=m(qg(\mathbf{x}))=(mq)g(\mathbf{x}),
\end{equation*}
i.e. $mq\in M^{\mathbf{x}}.$ If $M\in \mathcal{M}^{\mathcal{C}},$ then we
have for all $m\in M$ and $q\in Q:$%
\begin{equation*}
\begin{tabular}{lllll}
$\varrho _{M}(mq)$ & $=$ & $\varrho _{M}(\sum m_{<0>}q(m_{<1>}))$ &  &  \\ 
& $=$ & $\sum m_{<0><0>}\otimes _{A}m_{<0><1>}q(m_{<1>})$ &  &  \\ 
& $=$ & $\sum m_{<0>}\otimes _{A}m_{<1>1}q(m_{<1>2})$ &  &  \\ 
& $=$ & $\sum m_{<0>}\otimes _{A}q(m_{<1>})\mathbf{x}$ &  &  \\ 
& $=$ & $\sum m_{<0>}q(m_{<1>})\otimes _{A}\mathbf{x}$ &  &  \\ 
& $=$ & $\sum mq\otimes _{A}\mathbf{x},$ &  & 
\end{tabular}
\end{equation*}
i.e. $mq\in M^{co\mathcal{C}}.\blacksquare $
\end{enumerate}
\end{Beweis}

\begin{lemma}
\label{bim}

\begin{enumerate}
\item  With the canonical actions $A$ is a $(B,^{\ast }\mathcal{C})$%
-bimodule.

\item  $Q$ is a $(^{\ast }\mathcal{C},B)$-bimodule.
\end{enumerate}
\end{lemma}

\begin{Beweis}
\begin{enumerate}
\item  By assumption $A\in \mathcal{M}^{\mathcal{C}}\subseteq \mathcal{M}%
_{^{\ast }\mathcal{C}}.$ For all $b\in B,$ $a\in A$ and $g\in $ $^{\ast }%
\mathcal{C}$ we have 
\begin{equation*}
b(a\leftharpoonup g)=bg(\mathbf{x}a)=g(b(\mathbf{x}a))=g(\mathbf{x}%
(ba))=(ba)\leftharpoonup g.
\end{equation*}

\item  For all $a\in A,$ $q\in Q$ and $c\in \mathcal{C}$ we have 
\begin{equation*}
\sum c_{1}(aq)(c_{2})=\sum c_{1}q(c_{2}a)=\sum (ca)_{1}q((ca)_{2})=q(ca)%
\mathbf{x}=(aq)(c)\mathbf{x}.
\end{equation*}

For all $q\in Q,$ $b\in B$ and $c\in \mathcal{C}$ we have 
\begin{equation*}
\sum c_{1}(qb)(c_{2})=\sum c_{1}q(c_{2})b=q(c)\mathbf{x}b=q(c)b\mathbf{x}%
=(qb)(c)\mathbf{x}.
\end{equation*}
On the other hand we have for all $q\in Q,$ $g\in $ $^{\ast }\mathcal{C}$
and $c\in \mathcal{C}:$%
\begin{equation*}
\begin{tabular}{lllll}
$\sum c_{1}(g\cdot q)(c_{2})$ & $=$ & $\sum c_{1}q(c_{21}g(c_{22}))$ & $=$ & 
$\sum c_{11}q(c_{12}g(c_{2}))$ \\ 
& $=$ & $\sum c_{11}(g(c_{2})q)(c_{12})$ & $=$ & $\sum (g(c_{2})q)(c_{1})%
\mathbf{x}$ \\ 
& $=$ & $\sum q(c_{1}g(c_{2}))\mathbf{x}$ & $=$ & $(g\cdot q)(c)\mathbf{x}.$%
\end{tabular}
\end{equation*}

Moreover we have for all $b\in B,$ $q\in Q,$ $g\in $ $^{\ast }\mathcal{C}$
and $c\in \mathcal{C}:$%
\begin{equation*}
\begin{tabular}{lllll}
$((g\cdot q)b)(c)$ & $=$ & $(g\cdot q)(c)b$ & $=$ & $\sum q(c_{1}g(c_{2}))b$
\\ 
& $=$ & $\sum (qb)(c_{1}g(c_{2}))$ & $=$ & $(g\cdot qb)(c).\blacksquare $%
\end{tabular}
\end{equation*}
\end{enumerate}
\end{Beweis}

\begin{theorem}
\label{Mc}Keep the notation above fixed.

\begin{enumerate}
\item  $(A^{\mathbf{x}},^{\ast }\mathcal{C},A,(^{\ast }\mathcal{C})^{\mathbf{%
x}},\widetilde{F},\widetilde{G})$ is a Morita context derived form $%
A_{^{\ast }\mathcal{C}},$ where 
\begin{equation*}
\begin{tabular}{llllllll}
$\widetilde{F}$ & $:$ & $(^{\ast }\mathcal{C})^{\mathbf{x}}\otimes _{A^{%
\mathbf{x}}}A$ & $\longrightarrow $ & $^{\ast }\mathcal{C},$ & $q\otimes
_{A^{\mathbf{x}}}a$ & $\mapsto $ & $qa,$ \\ 
$\widetilde{G}$ & $:$ & $A\otimes _{^{\ast }\mathcal{C}}(^{\ast }\mathcal{C}%
)^{\mathbf{x}}$ & $\longrightarrow $ & $A^{\mathbf{x}},$ & $a\otimes
_{^{\ast }\mathcal{C}}q$ & $\mapsto $ & $a\leftharpoonup q.$%
\end{tabular}
\end{equation*}

\item  $(B,^{\ast }\mathcal{C},A,Q,F,G)$ is a Morita context, where 
\begin{equation*}
\begin{tabular}{llllllll}
$F$ & $:$ & $Q\otimes _{B}A$ & $\longrightarrow $ & $^{\ast }\mathcal{C},$ & 
$q\otimes _{B}a$ & $\mapsto $ & $qa,$ \\ 
$G$ & $:$ & $A\otimes _{^{\ast }\mathcal{C}}Q$ & $\longrightarrow $ & $B,$ & 
$a\otimes _{^{\ast }\mathcal{C}}q$ & $\mapsto $ & $a\leftharpoonup q.$%
\end{tabular}
\end{equation*}

If moreover $_{A}\mathcal{C}$ is locally projective, then the two Morita
contexts coincide.
\end{enumerate}
\end{theorem}

\begin{Beweis}
\begin{enumerate}
\item  By Lemma \ref{Psi-pr} we have $\mathrm{End}(A_{^{\ast }\mathcal{C}%
})\simeq A^{\mathbf{x}},$ $(^{\ast }\mathcal{C})^{\mathbf{x}}\simeq \mathrm{%
Hom}_{-^{\ast }\mathcal{C}}(A,^{\ast }\mathcal{C})$ and the result follows
by \cite[Proposition 12.6]{Fai81}.

\item  By Lemma \ref{bim} $A$ is a $(B,^{\ast }\mathcal{C})$-bimodule and $Q$
is a $(^{\ast }\mathcal{C},B)$-bimodule. For all $q\in Q,g\in $ $^{\ast }%
\mathcal{C},$ $a\in A$ and $c\in \mathcal{C}$ we have 
\begin{equation*}
F(g\cdot q\otimes _{B}a)(c)=\sum q(c_{2}g(c_{1}))a=(g\cdot qa)(c)=(g\cdot
F(q\otimes _{B}a))(c)
\end{equation*}
and 
\begin{equation*}
\begin{tabular}{lllll}
$F(q\otimes _{B}a\leftharpoonup g)(c)$ & $=$ & $q(c)(a\leftharpoonup g)$ & $%
= $ & $q(c)g(\mathbf{x}a)$ \\ 
& $=$ & $g(q(c)\mathbf{x}a)$ & $=$ & $\sum g(c_{1}q(c_{2})a)$ \\ 
& $=$ & $\sum g(c_{1}(qa)(c_{2}))$ & $=$ & $(F(q\otimes _{B}a)\cdot g)(c),$%
\end{tabular}
\end{equation*}
\newline
hence $F$ is $^{\ast }\mathcal{C}$-bilinear. Note that by Lemma \ref{Psi-pr} 
$G$ is well defined and is obviously $B$-bilinear. Moreover we have for all $%
a,\widetilde{a}\in A$ and $q,\widetilde{q}\in Q$ the following associativity
relations: 
\begin{equation*}
\begin{tabular}{lllll}
$(F(q\otimes _{B}a)\cdot \widetilde{q})(c)$ & $=$ & $\sum \widetilde{q}%
(c_{1}q(c_{2})a)$ & $=$ & $\widetilde{q}(q(c)\mathbf{x}a)$ \\ 
& $=$ & $q(c)\widetilde{q}(\mathbf{x}a)$ & $=$ & $(qG(a\otimes _{^{\ast }%
\mathcal{C}}\widetilde{q}))(c),$ \\ 
$G(a\otimes _{^{\ast }\mathcal{C}}q)\widetilde{a}$ & $=$ & $q(\mathbf{x}a)%
\widetilde{a}$ & $=$ & $(q\widetilde{a})(\mathbf{x}a)$ \\ 
& $=$ & $F(q\otimes _{B}\widetilde{a})(\mathbf{x}a)$ & $=$ & $%
a\leftharpoonup F(q\otimes _{B}\widetilde{a}).$%
\end{tabular}
\end{equation*}
\newline
If $_{A}\mathcal{C}$ is locally projective, then $A^{\mathbf{x}}=A^{co%
\mathcal{C}},$ $(^{\ast }\mathcal{C})^{\mathbf{x}}=Q$ by Lemma \ref{Psi-pr}
and the two contexts coincide.$\blacksquare $
\end{enumerate}
\end{Beweis}

\begin{punto}
\cite[Definition 5.3]{Brz02} An $A$-coring $\mathcal{C}$ is said to be 
\textbf{Galois}, if there exists an $A$-coring isomorphism $\chi :A\otimes
_{B}A\longrightarrow \mathcal{C}$ such that $\chi (1_{A}\otimes _{B}1_{A})=%
\mathbf{x}.$ Recall that $A\otimes _{B}A$ is an $A$-coring with the
canonical $A$-bimodule structure, comultiplication 
\begin{equation*}
\Delta :A\otimes _{B}A\longrightarrow (A\otimes _{B}A)\otimes _{A}(A\otimes
_{B}A),\text{ }\widetilde{a}\otimes _{B}a\mapsto (\widetilde{a}\otimes
_{B}1_{A})\otimes _{A}(1_{A}\otimes _{B}a)
\end{equation*}
and counity $\varepsilon _{A\otimes _{B}A}:A\otimes _{B}A\longrightarrow A,$ 
$\widetilde{a}\otimes _{B}a\mapsto \widetilde{a}a.$
\end{punto}

\begin{punto}
Consider the functors 
\begin{equation*}
(-)^{co\mathcal{C}}:\mathcal{M}^{\mathcal{C}}\longrightarrow \mathcal{M}_{B}%
\text{ and }-\otimes _{B}A:\mathcal{M}_{B}\longrightarrow \mathcal{M}^{%
\mathcal{C}}.
\end{equation*}
By \cite[Proposition 5.2]{Brz02} $(-\otimes _{B}A,(-)^{co\mathcal{C}})$ is
an adjoint pair of covariant functors, where the adjunctions are given by 
\begin{equation}
\Phi _{N}:N\longrightarrow (N\otimes _{B}A)^{co\mathcal{C}},\text{ }n\mapsto
n\otimes _{B}1_{A}  \label{str}
\end{equation}
and 
\begin{equation}
\Psi _{M}:M^{co\mathcal{C}}\otimes _{B}A\longrightarrow M,\text{ }m\otimes
_{B}a\mapsto ma.  \label{wk}
\end{equation}
If $\Psi _{M}$ is an isomorphism for all $M\in \mathcal{M}^{\mathcal{C}},$
then we say $\mathcal{M}^{\mathcal{C}}$ satisfies the \textbf{weak structure
theorem. }If in addition $\Phi _{N}$ is an isomorphism for all $N\in 
\mathcal{M}_{B},$ then we say $\mathcal{M}^{\mathcal{C}}$ satisfies the 
\textbf{strong structure theorem }(in this case $(-)^{co\mathcal{C}}$ and $%
-\otimes _{B}A$ give an equivalence of categories $\mathcal{M}^{\mathcal{C}%
}\simeq \mathcal{M}_{B}$).
\end{punto}

\begin{punto}
Let $W\in \mathcal{M}_{A}$ and consider the canonical right $\mathcal{C}$%
-comodule $W\otimes _{A}\mathcal{C}.$ Then $W\simeq (W\otimes _{A}\mathcal{C}%
)^{co\mathcal{C}}$ via $w\mapsto w\otimes _{A}\mathbf{x}$ with inverse $%
w\otimes _{A}c\mapsto w\varepsilon _{\mathcal{C}}(c)$ and we define 
\begin{equation}
\beta _{W}:=\Psi _{W\otimes _{A}\mathcal{C}}:W\otimes _{B}A\longrightarrow
W\otimes _{A}\mathcal{C},\text{ }w\otimes _{B}a\mapsto w\otimes _{A}\mathbf{x%
}a.  \label{bw}
\end{equation}
In particular we have for $W=A$ the morphism of $A$-corings 
\begin{equation}
\beta :=\Psi _{A\otimes _{A}\mathcal{C}}:A\otimes _{B}A\longrightarrow
A\otimes _{A}\mathcal{C}\simeq \mathcal{C},\text{ }\widetilde{a}\otimes
_{B}a\mapsto \widetilde{a}\mathbf{x}a.  \label{bet}
\end{equation}
If $\beta $ is bijective, then $\mathcal{C}$ is a Galois $A$-coring and we
call the ring extension $A/B$ $\mathcal{C}$-\textbf{Galois}.
\end{punto}

\begin{theorem}
\label{surj}For the Morita context $(B,^{\ast }\mathcal{C},A,Q,F,G)$ the
following statements are equivalent:

\begin{enumerate}
\item  $G:A\otimes _{^{\ast }\mathcal{C}}Q\longrightarrow B$ is surjective 
\emph{(}bijective and $B=A^{\mathbf{x}}$\emph{)};

\item  there exists $\widehat{q}\in Q,$ such that $\widehat{q}(\mathbf{x}%
)=1_{A};$

\item  for every right $^{\ast }\mathcal{C}$-module $M$ we have a $B$-module
isomorphism $M\otimes _{^{\ast }\mathcal{C}}Q\simeq M^{\mathbf{x}}.$

\item  for every right $\mathcal{C}$-comodule $M$ we have $M\otimes _{^{\ast
}\mathcal{C}}Q\simeq M^{\mathbf{co}\mathcal{C}}$ as $B$-modules.

If moreover $_{A}\mathcal{C}$ is locally projective, then \emph{(1)-(4)}\
are moreover equivalent to:

\item  $A_{^{\ast }\mathcal{C}}$ is \emph{(}f.g.\emph{)} projective.
\end{enumerate}
\end{theorem}

\begin{Beweis}
(1) $\Rightarrow $ (2). Assume $G$ to be surjective. Then there exist $%
a_{1},...,a_{k}$ and $q_{1},...,q_{k}\in Q,$ such that $G(\sum%
\limits_{i=1}^{k}a_{i}\otimes _{^{\ast }\mathcal{C}}q_{i})=1_{A}.$ Set $%
\widehat{q}:=\sum\limits_{i=1}^{k}a_{i}q_{i}\in Q.$ Then we have 
\begin{equation*}
\widehat{q}(\mathbf{x})=(\sum\limits_{i=1}^{k}a_{i}q_{i})(\mathbf{x}%
)=\sum\limits_{i=1}^{k}q_{i}(\mathbf{x}a_{i})=\sum\limits_{i=1}^{k}(a_{i}%
\leftharpoonup q_{i})=G(\sum\limits_{i=1}^{k}a_{i}\otimes _{^{\ast }\mathcal{%
C}}q_{i})=1_{A}.
\end{equation*}

(2)\ $\Rightarrow $ (3). Consider the $B$-module morphism 
\begin{equation*}
\xi _{M}:M\otimes _{^{\ast }\mathcal{C}}Q\longrightarrow M^{\mathbf{x}},%
\text{ }m\otimes _{^{\ast }\mathcal{C}}q\mapsto mq.
\end{equation*}
Let $\widehat{q}\in Q$ with $\widehat{q}(\mathbf{x})=1_{A}$ and define $%
\widetilde{\xi }_{M}:M^{\mathbf{x}}\longrightarrow M\otimes _{^{\ast }%
\mathcal{C}}Q,$ $m\mapsto m\otimes _{^{\ast }\mathcal{C}}\widehat{q}.$ For
every $n\in M^{\mathbf{x}}$ we have 
\begin{equation*}
(\xi _{M}\circ \widetilde{\xi }_{M})(n)=\xi _{M}(n\otimes _{^{\ast }\mathcal{%
C}}\widehat{q})=n\leftharpoonup \widehat{q}=n\widehat{q}(\mathbf{x})=n.
\end{equation*}
On the other hand we have for all $m\in M$ and $q\in Q:$ 
\begin{equation*}
\begin{tabular}{lllllll}
$(\widetilde{\xi }_{M}\circ \xi _{M})(m\otimes _{^{\ast }\mathcal{C}}q)$ & $%
= $ & $\widetilde{\xi }_{M}(m\leftharpoonup q)$ & $=$ & $m\leftharpoonup
q\otimes _{^{\ast }\mathcal{C}}\widehat{q}$ & $=$ & $m\otimes _{^{\ast }%
\mathcal{C}}q\cdot \widehat{q}$ \\ 
& $=$ & $m\otimes _{^{\ast }\mathcal{C}}q\widehat{q}(\mathbf{x})$ & $=$ & $%
m\otimes _{^{\ast }\mathcal{C}}q,$ &  & 
\end{tabular}
\end{equation*}
i.e. $\xi _{M}$ is bijective with inverse $\widetilde{\xi }_{M}.$

(3)\ $\Rightarrow $ (4). Let $M\in \mathcal{M}^{\mathcal{C}}.$ By Lemma \ref
{Psi-pr} we have $\xi _{M}(M\otimes _{^{\ast }\mathcal{C}}Q)\subseteq M^{co%
\mathcal{C}}\subseteq M^{\mathbf{x}}.$ By assumption $\xi _{M}:A\otimes
_{^{\ast }\mathcal{C}}Q\longrightarrow M^{\mathbf{x}}$ is bijective. Hence $%
M^{\mathbf{x}}=M^{co\mathcal{C}}$ and $M\otimes _{^{\ast }\mathcal{C}}Q%
\overset{\xi _{M}}{\simeq }M^{co\mathcal{C}}.$

(4)\ $\Rightarrow $ (1) We are done since $G=\xi _{A}.$

Assume $_{A}\mathcal{C}$ to be locally projective.

Then $B\simeq \mathrm{End}(A_{^{\ast }\mathcal{C}}),$ $Q\simeq \mathrm{Hom}%
_{-^{\ast }\mathcal{C}}(A,^{\ast }\mathcal{C})$ and we get (1) $%
\Longleftrightarrow $ (5) by \cite[Corollary 12.8]{Fai81}.$\blacksquare $
\end{Beweis}

\begin{corollary}
\label{co=x}For the Morita context $(B,^{\ast }\mathcal{C},A,Q,F,G)$ assume
there exists $\widehat{q}\in Q$ with $\widehat{q}(\mathbf{x})=1_{A}$ \emph{(}%
equivalently $G:Q\otimes _{B}A\longrightarrow $ $^{\ast }\mathcal{C}$ is
surjective\emph{)}. Then:

\begin{enumerate}
\item  For every $N\in \mathcal{M}_{B},$ $\Phi _{N}$ is an isomorphism.

\item  $B$ is a left $B$-direct summand of $A.$
\end{enumerate}
\end{corollary}

\begin{Beweis}
\begin{enumerate}
\item  Let $N\in \mathcal{M}_{B}.$ Then we have by Theorem \ref{surj} the
isomorphisms $G:A\otimes _{^{\ast }\mathcal{C}}Q\longrightarrow B$ and $\xi
_{N\otimes _{B}A}:(N\otimes _{B}A)\otimes _{^{\ast }\mathcal{C}%
}Q\longrightarrow (N\otimes _{B}A)^{co\mathcal{C}}.$ Moreover $\Phi _{N}$ is
given by the canonical isomorphisms 
\begin{equation*}
N\simeq N\otimes _{B}B\simeq N\otimes _{B}(A\otimes _{^{\ast }\mathcal{C}%
}Q)\simeq (N\otimes _{B}A)\otimes _{^{\ast }\mathcal{C}}Q\simeq (N\otimes
_{B}A)^{co\mathcal{C}}.
\end{equation*}

\item  The map $\mathrm{tr}_{A}:A\longrightarrow B,$ $a\mapsto
a\leftharpoonup \widehat{q}$ is left $B$-linear with $\mathrm{tr}_{A}(b)=b$
for all $b\in B.\blacksquare $
\end{enumerate}
\end{Beweis}

\begin{corollary}
\label{general}For the Morita context $(B,^{\ast }\mathcal{C},A,Q,F,G)$
assume there exists $\widehat{q}\in Q$ with $\widehat{q}(\mathbf{x})=1_{A}$ 
\emph{(}equivalently $G:Q\otimes _{B}A\longrightarrow B$ is surjective\emph{)%
}. Then:

\begin{enumerate}
\item  $_{B}A$ and $Q_{B}$ are generators.

\item  $A_{^{\ast }\mathcal{C}}$ and $_{^{\ast }\mathcal{C}}Q$ are f.g. and
projective.

\item  $F:Q\otimes _{B}A\longrightarrow $ $^{\ast }\mathcal{C}$ induces
bimodule isomorphisms 
\begin{equation*}
A\simeq \mathrm{Hom}_{^{\ast }\mathcal{C}-}(Q,^{\ast }\mathcal{C})\text{ and 
}Q\simeq \mathrm{Hom}_{-^{\ast }\mathcal{C}}(A,^{\ast }\mathcal{C}).
\end{equation*}

\item  The bimodule structures above induce ring isomorphisms 
\begin{equation*}
B\simeq \mathrm{End}(A_{^{\ast }\mathcal{C}})\text{ and }B\simeq \mathrm{End}%
(_{^{\ast }\mathcal{C}}Q)^{op}.
\end{equation*}
\end{enumerate}
\end{corollary}

\begin{Beweis}
The result follows by standard argument of Morita Theory (e.g. 
\cite[Proposition 12.7]{Fai81}).$\blacksquare $
\end{Beweis}

\begin{proposition}
\label{f-G}Consider the Morita context $(B,^{\ast }\mathcal{C},A,Q,F,G)$ and
assume that $F:Q\otimes _{^{\ast }\mathcal{C}}A\longrightarrow $ $^{\ast }%
\mathcal{C}$ is surjective. Then:

\begin{enumerate}
\item  $A_{^{\ast }\mathcal{C}}$ is a generator, $Q\simeq \mathrm{Hom}%
_{B-}(A,B)$ as bimodules and $^{\ast }\mathcal{C}\simeq \mathrm{End}(Q_{B}).$

\item  $\mathcal{M}^{\mathcal{C}}$ satisfies the weak structure theorem 
\emph{(}in particular $A/B$ is $\mathcal{C}$-Galois\emph{)}.
\end{enumerate}
\end{proposition}

\begin{Beweis}
\begin{enumerate}
\item  The result$\;$follows by standard argument of Morita Theory (e.g. 
\cite[Proposition 12.7]{Fai81}).

\item  By assumption $\varepsilon _{\mathcal{C}}=F(\sum%
\limits_{i=1}^{k}q_{i}\otimes _{B}a_{i})$ for some $\{(q_{i},a_{i})%
\}_{i=1}^{k}\subseteq Q\times A.$ In this case $\Psi _{M}:M^{co\mathcal{C}%
}\otimes _{B}A\longrightarrow M$ is bijective with inverse $\widetilde{\Psi }%
_{M}:M\longrightarrow M^{co\mathcal{C}}\otimes _{B}A,$ $m\mapsto
\sum\limits_{i=1}^{k}mq_{i}\otimes _{B}a_{i}.$ In fact, we have for all $%
m\in M,$ $n\in M^{co\mathcal{C}}$ and $a\in A:$%
\begin{equation*}
\begin{tabular}{lllll}
$(\Psi _{M}\circ \widetilde{\Psi }_{M})(m)$ & $=$ & $\sum%
\limits_{i=1}^{k}(mq_{i})a_{i}$ & $=$ & $\sum%
\limits_{i=1}^{k}(m_{<0>}q_{i}(m_{<1>})a_{i}$ \\ 
& $=$ & $\sum\limits_{i=1}^{k}m_{<0>}(q_{i}a_{i})(m_{<1>})$ & $=$ & $%
\sum\limits_{i=1}^{k}m_{<0>}\varepsilon _{\mathcal{C}}(m_{<1>})$ \\ 
& $=$ & $m$ &  & 
\end{tabular}
\end{equation*}
and 
\begin{equation*}
\begin{tabular}{lllll}
$(\widetilde{\Psi }_{M}\circ \Psi _{M})(n\otimes _{B}a)$ & $=$ & $%
\sum\limits_{i=1}^{k}(na)q_{i}\otimes _{B}a_{i}$ & $=$ & $%
\sum\limits_{i=1}^{k}nq_{i}(\mathbf{x}a)\otimes _{B}a_{i}$ \\ 
& $=$ & $\sum\limits_{i=1}^{k}n\otimes _{B}q_{i}(\mathbf{x}a)a_{i}$ & $=$ & $%
\sum\limits_{i=1}^{k}n\otimes _{B}(q_{i}a_{i})(\mathbf{x}a)$ \\ 
& $=$ & $n\otimes _{B}\varepsilon _{\mathcal{C}}(\mathbf{x}a)$ & $=$ & $%
n\otimes _{B}a.\blacksquare $%
\end{tabular}
\end{equation*}
\end{enumerate}
\end{Beweis}

\begin{theorem}
\label{C -finite}For the Morita context $(B,^{\ast }\mathcal{C},A,Q,F,G)$
the following are equivalent:

\begin{enumerate}
\item  $F:Q\otimes _{B}A\longrightarrow $ $^{\ast }\mathcal{C}$ is
surjective \emph{(}bijective\emph{)};

\item  (a)\ $Q_{B}$ is f.g. and projective;

(b) $\Omega :A\longrightarrow \mathrm{Hom}_{-B}(Q,B),$ $a\mapsto \lbrack
q\mapsto a\leftharpoonup q]$ is a bimodule isomorphism;

(c) $_{^{\ast }\mathcal{C}}Q$ is faithful.

If $_{A}\mathcal{C}$ is $A$-cogenerated, then \emph{(1)} $\&$ \emph{(2)} are
moreover equivalent to:

\item  (a)$\;_{B}A$ is f.g. and projective;

(b)\ $\Lambda :$ $^{\ast }\mathcal{C}\longrightarrow \mathrm{End}%
(_{B}A)^{op},$ $g\mapsto \lbrack a\mapsto a\leftharpoonup g]$ is a ring
isomorphism.

\item  $A_{^{\ast }\mathcal{C}}$ is a generator.

If moreover $_{A}\mathcal{C}$ is f.g. and projective, then \emph{(1)-(4)}
are equivalent to:

\item  $\mathcal{M}^{\mathcal{C}}$ satisfies the weak structure theorem.
\end{enumerate}
\end{theorem}

\begin{Beweis}
The implications (1)\ $\Rightarrow $ (2), (3), (4) follow without any
finiteness conditions on $\mathcal{C}$ by standard argument of Morita Theory
(e.g. \cite[Proposition 12.7]{Fai81}). Note that $_{^{\ast }\mathcal{C}}Q$
is faithful by the embedding $^{\ast }\mathcal{C}\hookrightarrow \mathrm{End}%
(Q_{B})$ (see Proposition \ref{f-G} (1)).

(2)\ $\Rightarrow $ (1). Let $\{(q_{i},p_{i})\}_{i=1}^{k}\subset Q\times 
\mathrm{Hom}_{-B}(Q,B)$ be a dual basis for $Q_{B}.$ By (b)\ there exist $%
a_{1},...,a_{k}\in A,$ such that $\Omega (a_{i})=q_{i}$ for $i=1,...,k.$ For
every $q\in Q$ we have then $(\sum\limits_{i=1}^{k}q_{i}a_{i})\cdot
q=\sum\limits_{i=1}^{k}q_{i}(a_{i}\leftharpoonup
q)=\sum\limits_{i=1}^{k}q_{i}p_{i}(q)=q,$ hence $\sum%
\limits_{i=1}^{k}q_{i}a_{i}=\varepsilon _{\mathcal{C}}$ by (c) and the $%
^{\ast }\mathcal{C}$-bilinear morphism $F:Q\otimes _{^{\ast }\mathcal{C}%
}A\longrightarrow $ $^{\ast }\mathcal{C}$ is surjective.

Assume $_{A}\mathcal{C}$ to be $A$-cogenerated.

(3) $\Rightarrow $ (1). Let $\{(a_{i},p_{i})\}_{i=1}^{k}\subset A\times 
\mathrm{Hom}_{B-}(A,B)$ be a dual basis of $_{B}A.$ By (b), there exist $%
g_{1},...,g_{k}\in $ $^{\ast }\mathcal{C},$ such that $\Lambda (g_{i})=p_{i}$
for $i=1,...,k.$ \textbf{Claim}: $g_{1},...,g_{k}\in Q.$ For all $f\in $ $%
^{\ast }\mathcal{C}$ and $i=1,...,k$ we have 
\begin{equation*}
\begin{tabular}{lllll}
$\Lambda (g_{i}\cdot f)(a)$ & $=$ & $a\leftharpoonup (g_{i}\cdot f)$ & $=$ & 
$(a\leftharpoonup g_{i})\leftharpoonup f)$ \\ 
& $=$ & $p_{i}(a)\leftharpoonup f$ & $=$ & $f(\mathbf{x}p_{i}(a))$ \\ 
& $=$ & $f(p_{i}(a)\mathbf{x})$ & $=$ & $p_{i}(a)f(\mathbf{x})$ \\ 
& $=$ & $(p_{i}f(\mathbf{x}))(a)$ & $=$ & $\Lambda (g_{i}f(\mathbf{x}))(a),$%
\end{tabular}
\end{equation*}
hence $g_{i}\cdot f=g_{i}f(\mathbf{x}),$ i.e. $g_{i}\in $ $(^{\ast }\mathcal{%
C})^{\mathbf{x}}=Q$ (by Lemma \ref{Psi-pr} (2)). Moreover for every $a\in A$
we have: $\Lambda
(\sum\limits_{i=1}^{k}g_{i}a_{i})(a)=\sum\limits_{i=1}^{k}a\leftharpoonup
g_{i}a_{i}=\sum\limits_{i=1}^{k}p_{i}(a)a_{i}=a,$ i.e. $\sum%
\limits_{i=1}^{k}g_{i}a_{i}=\varepsilon _{\mathcal{C}}$ and the $^{\ast }%
\mathcal{C}$-bilinear morphism $F$ is surjective.

(4)\ $\Rightarrow $ (1). Since $Q\simeq \mathrm{Hom}_{-\text{ }^{\ast }%
\mathcal{C}}(A,^{\ast }\mathcal{C}),$ we have $\func{Im}(F)=\mathrm{tr}%
(A,^{\ast }\mathcal{C}):=\sum \{\func{Im}(h):h\in \mathrm{Hom}_{-\text{ }%
^{\ast }\mathcal{C}}(A,^{\ast }\mathcal{C})\},$ hence $\func{Im}(F)=$ $%
^{\ast }\mathcal{C}$ iff $A_{^{\ast }\mathcal{C}}$ is a generator (e.g. 
\cite[Page 154]{Wis88}).

Assume $_{A}\mathcal{C}$ to be f.g. and projective.

(1)\ $\Rightarrow $ (5)\ follows without any finiteness conditions on $%
\mathcal{C}$ by Proposition \ref{f-G} (2).

(5) $\Rightarrow $ (1). Since $_{A}\mathcal{C}$ is f.g. and projective, we
have $\mathcal{M}^{\mathcal{C}}\simeq \mathcal{M}_{^{\ast }\mathcal{C}}$
(e.g. \cite[Lemma 4.3]{Brz02}), hence $^{\ast }\mathcal{C}\in \mathcal{M}^{%
\mathcal{C}},$ $Q=(^{\ast }\mathcal{C})^{co\mathcal{C}}$ and $F=\Psi
_{^{\ast }\mathcal{C}}.\blacksquare $
\end{Beweis}

\section{Galois Extensions and Equivalences}

The notation of the first section remains fixed. For every $M\in \mathcal{M}%
^{\mathcal{C}}$ we have the $\mathcal{C}$-colinear morphism 
\begin{equation*}
\Psi _{M}^{\prime }:\mathrm{Hom}^{\mathcal{C}}(A,M)\otimes
_{B}A\longrightarrow M,\text{ }f\otimes _{B}a\mapsto f(a).
\end{equation*}

In this section we characterize $A$ being a generator (resp. a progenerator)
in $\mathcal{M}^{\mathcal{C}}$ under the assumption that $_{A}\mathcal{C}$
is locally projective. Our approach is similar to that of \cite{MZ97} and
our results generalize those obtained there for the special case of the
category of Doi-Koppinen modules $\mathcal{M}(H)_{A}^{C}.$

\begin{lemma}
\label{A-Gen}Assume $_{A}\mathcal{C}$ to be locally projective. If $_{B}A$
is flat and $A/B$ is $\mathcal{C}$-Galois, then:

\begin{enumerate}
\item  $A$ is a subgenerator in $\mathcal{M}^{\mathcal{C}},$ i.e. $\sigma
\lbrack A_{^{\ast }\mathcal{C}}]=\sigma \lbrack \mathcal{C}_{^{\ast }%
\mathcal{C}}].$

\item  for each $M\in \mathcal{M}^{\mathcal{C}},$ $\Psi _{M}^{\prime }$ is
injective.

\item  for every $A$-generated $M\in \mathcal{M}^{\mathcal{C}},$ $\Psi
_{M}^{\prime }$ is an isomorphism.
\end{enumerate}
\end{lemma}

\begin{Beweis}
Assume $_{A}\mathcal{C}$ to be locally projective.

\begin{enumerate}
\item  Since $A/B$ is $\mathcal{C}$-Galois, $\beta ^{\prime }:=\Psi _{%
\mathcal{C}}^{\prime }$ is an isomorphism, hence $\mathcal{C}$ is $A$%
-generated. Consequently $\sigma \lbrack A_{^{\ast }\mathcal{C}}]\subseteq
\sigma \lbrack \mathcal{C}_{^{\ast }\mathcal{C}}]\subseteq \sigma \lbrack
A_{^{\ast }\mathcal{C}}],$ i.e. $\sigma \lbrack A_{^{\ast }\mathcal{C}%
}]=\sigma \lbrack \mathcal{C}_{^{\ast }\mathcal{C}}].$

\item  With slight modifications, the proof of \cite[Lemma 3.22]{MZ97}
applies.

\item  If $M\in \mathcal{M}^{\mathcal{C}}$ is $A$-generated, then $\Psi
_{M}^{\prime }$ is surjective, hence bijective by (2).$\blacksquare $
\end{enumerate}
\end{Beweis}

\qquad The following result is a generalization of \cite[Proposition 3.13]
{Brz99} (which in turn generalizes \cite[Theorem 2.11]{DT89}):

\begin{proposition}
\label{suf-con}Assume $A/B$ to be $\mathcal{C}$-Galois.

\begin{enumerate}
\item  If $_{B}A$ is flat, then $\mathcal{M}^{\mathcal{C}}$ satisfies the
weak structure theorem.

\item  Assume there exists $\widehat{q}\in Q,$ such that $\widehat{q}(%
\mathbf{x})=1_{A}.$ If $_{B}A$ is flat, or for all $b\in B$ and $c\in 
\mathcal{C}$ we have $\widehat{q}(cb)=q(c)b,$ then $\mathcal{M}^{\mathcal{C}%
} $ satisfies the strong structure theorem.
\end{enumerate}
\end{proposition}

\begin{Beweis}
\begin{enumerate}
\item  The proof is the first part of the proof of \cite[Theorem 5.6]{Brz02}.

\item  By assumption and Corollary \ref{co=x}, $\Phi _{N}$ is an isomorphism
for all $N\in \mathcal{M}_{B}.$ If $_{B}A$ is flat, then $\mathcal{M}^{%
\mathcal{C}}$ satisfies the weak structure theorem by (1). On the other
hand, if for all $b\in B$ and $c\in \mathcal{C}$ we have $\widehat{q}%
(cb)=q(c)b,$ then an analog argument to that in the proof of 
\cite[Proposition 3.13]{Brz99} shows that $\mathcal{M}^{\mathcal{C}}$
satisfies the weak structure theorem.$\blacksquare $
\end{enumerate}
\end{Beweis}

\begin{theorem}
\label{gen}Assume $_{A}\mathcal{C}$ to be locally projective. Then the
following are equivalent:

\begin{enumerate}
\item  $\mathcal{M}^{\mathcal{C}}$ satisfies the weak structure theorem;

\item  $_{B}A$ is flat and $A/B$ is $\mathcal{C}$-Galois;

\item  $_{B}A$ is flat and $\beta ^{\prime }:=\Psi _{\mathcal{C}}^{\prime }$
is an isomorphism;

\item  $_{B}A$ is flat and for every $A$-generated $M\in \mathcal{M}^{%
\mathcal{C}},$ $\Psi _{M}^{\prime }$ is bijective;

\item  for every $M\in \mathcal{M}^{\mathcal{C}}=\sigma \lbrack \mathcal{C}%
_{^{\ast }\mathcal{C}}],$ the $\mathcal{C}$-colinear morphism $\Psi
_{M}^{\prime }$ is bijective;

\item  $\sigma \lbrack \mathcal{C}_{^{\ast }\mathcal{C}}]=\mathrm{Gen}%
(A_{^{\ast }\mathcal{C}});$

\item  $_{B}A$ is flat, $\sigma \lbrack \mathcal{C}_{^{\ast }\mathcal{C}%
}]=\sigma \lbrack A_{^{\ast }\mathcal{C}}]$ and $\mathrm{Hom}_{-^{\ast }%
\mathcal{C}}(A,-):\mathrm{Gen}(A_{^{\ast }\mathcal{C}})\longrightarrow 
\mathcal{M}_{B}$ is full faithful;

\item  $\mathrm{Hom}^{\mathcal{C}}(A,-):\mathcal{M}^{\mathcal{C}%
}\longrightarrow \mathcal{M}_{B}$ is faithful;

\item  $A$ is a generator in $\mathcal{M}^{\mathcal{C}}.$
\end{enumerate}
\end{theorem}

\begin{Beweis}
(1) $\Longleftrightarrow $ (5)\ $\&\;$(2) $\Longleftrightarrow $ (3)\ follow
by Lemma \ref{Psi-pr}. The equivalences (4)\ $\Longleftrightarrow $ (5) $%
\Longleftrightarrow $ (6) $\Longleftrightarrow $ (7) follow by \cite[Theorem
2.3]{MZ97}. The equivalence (8) $\Longleftrightarrow $ (9)\ is evident for
any category, and moreover (6) $\Longleftrightarrow $ (9) by the fact that $%
\mathrm{Gen}(A_{^{\ast }\mathcal{C}})\subseteq \sigma \lbrack A_{^{\ast }%
\mathcal{C}}]\subseteq \sigma \lbrack \mathcal{C}_{^{\ast }\mathcal{C}}]=%
\mathcal{M}^{\mathcal{C}}.$ By Lemma \ref{A-Gen} we have (3) $\Rightarrow $
(4). Now assuming (1)\ we conclude that $A/B$ is $\mathcal{C}$-Galois and
that $_{B}A$ is flat (since (1) $\Longleftrightarrow $ (5) $%
\Longleftrightarrow $ (7)), hence (1)\ $\Rightarrow $ (2) follows and we are
done.$\blacksquare $
\end{Beweis}

\begin{definition}
(\cite[Definition 2.4]{MZ97})\ A left module $P$ over a ring $\mathcal{S}$
is called a \textbf{weak generator}, if for any right $\mathcal{S}$-module $%
Y,$ $Y\otimes _{\mathcal{S}}P=0$ implies $Y=0.$ A right module $P$ over a
ring $\mathcal{R}$ is called \textbf{quasiprogenerator} (resp. \textbf{%
progenerator}), if $P_{\mathcal{R}}$ is f.g. quasiprojective and generates
each of its submodules (resp. $P_{\mathcal{R}}$ is f.g., projective and a
generator). $P_{\mathcal{R}}$ is called \textbf{faithful }(resp. \textbf{%
balanced}), if the canonical morphism $\mathcal{R}\longrightarrow \mathrm{End%
}(_{\mathrm{End}(P_{\mathcal{R}})}P)^{op}$ is injective (resp. surjective).
\end{definition}

\begin{theorem}
\label{prog}Assume $_{A}\mathcal{C}$ to be flat. Then the following are
equivalent:

\begin{enumerate}
\item  $\mathcal{M}^{\mathcal{C}}$ satisfies the strong structure theorem;

\item  $_{B}A$ is faithfully flat and $A/B$ is $\mathcal{C}$-Galois.

If moreover $_{A}\mathcal{C}$ is locally projective, then \emph{(1)}\ $\&$ 
\emph{(2)} are moreover equivalent to:

\item  $_{B}A$ is faithfully flat and $\beta ^{\prime }:=\Psi _{\mathcal{C}%
}^{\prime }$ is bijective;

\item  $_{B}A$ is faithfully flat and for every $M\in \sigma \lbrack
A_{^{\ast }\mathcal{C}}],$ $\Psi _{M}^{\prime }$ is bijective;

\item  $A_{^{\ast }\mathcal{C}}$ is quasiprojective and generates each of
its submodules, $_{B}A$ is a weak generator and $\sigma \lbrack \mathcal{C}%
_{^{\ast }\mathcal{C}}]=\sigma \lbrack A_{^{\ast }\mathcal{C}}];$

\item  $A_{^{\ast }\mathcal{C}}$ is a quasiprogenerator and $\sigma \lbrack 
\mathcal{C}_{^{\ast }\mathcal{C}}]=\sigma \lbrack A_{^{\ast }\mathcal{C}}];$

\item  $_{B}A$ is a weak generator, $\Psi _{M}^{\prime }$ is an isomorphism
for every $M\in \mathrm{Gen}(A_{^{\ast }\mathcal{C}})$ and $\sigma \lbrack
A_{^{\ast }\mathcal{C}}]=\sigma \lbrack \mathcal{C}_{^{\ast }\mathcal{C}}];$

\item  $\mathrm{Hom}^{\mathcal{C}}(A,-):\mathcal{M}^{\mathcal{C}%
}\longrightarrow \mathcal{M}_{B}$ is an equivalence;

\item  $A$ is a progenerator in $\mathcal{M}^{\mathcal{C}}.$
\end{enumerate}
\end{theorem}

\begin{Beweis}
(1)\ $\Longleftrightarrow $(2)$\;$is \cite[Theorem 5.6]{Brz02}. Assume $_{A}%
\mathcal{C}$ to be locally projective. Then (2) $\Longleftrightarrow $ (3)\
follows by Lemma \ref{Psi-pr} and we get (1) $\Longleftrightarrow $ (8) $%
\Longleftrightarrow $ (9) by characterizations of progenerators in
categories of type $\sigma \lbrack M]$ (see \cite[18.5, 46.2]{Wis88}).
Moreover (4) $\Longleftrightarrow $ (5) $\Longleftrightarrow $ (6) $%
\Longleftrightarrow $ (7)\ follow from \cite[Theorem 2.5]{MZ97}. Obviously
(3) $\Rightarrow $ (4) (note that (3) $\Longleftrightarrow $ (2) $%
\Longleftrightarrow $ (1)). Assume now (4). Then $_{B}A$ is faithfully flat
and moreover $\Psi _{\mathcal{C}}^{\prime }$ is bijective, since $\mathcal{C}%
\in \sigma \lbrack A_{^{\ast }\mathcal{C}}]$ by (6), i.e. (4) $\Rightarrow $
(3)\ and the proof is complete.$\blacksquare $
\end{Beweis}

\begin{remark}
Assume $_{A}\mathcal{C}$ to be locally projective. Then $\func{Im}%
(F)\subseteq $ $^{\Box }\mathcal{C}.$ In fact we have for all $q\in Q,$ $%
a\in A,$ $g\in $ $^{\ast }\mathcal{C}$ and $c\in C:$ 
\begin{equation*}
((qa)\cdot g)(c)=\sum g(c_{1}q(c_{2})a)=g(q(c)\mathbf{x}a)=q(c)g(\mathbf{x}%
a)=(qg(\mathbf{x}a)(c),
\end{equation*}
hence $qa\in $ $^{\Box }\mathcal{C},$ with $\varrho (qa)=q\otimes _{A}%
\mathbf{x}a.$
\end{remark}

\begin{proposition}
Assume $_{A}\mathcal{C}$ to be locally projective and that there exists $%
\widehat{q}\in Q$ with $\widehat{q}(\mathbf{x})=1_{A}$ \emph{(}equivalently $%
G:A\otimes _{^{\ast }\mathcal{C}}Q\longrightarrow $ $B$ is surjective\emph{)}%
. Then $\mathcal{M}^{\mathcal{C}}$ satisfies the strong equivalence theorem,
iff $\func{Im}(F)=$ $^{\Box }\mathcal{C}$ and the following map is
surjective for every $M\in \mathcal{M}^{\mathcal{C}}$%
\begin{equation*}
\varpi _{M}:M\otimes _{^{\ast }\mathcal{C}}\text{ }^{\Box }\mathcal{C}%
\longrightarrow M,\text{ }m\otimes _{^{\ast }\mathcal{C}}f\mapsto mf.
\end{equation*}
In this case $Q\otimes _{B}A\overset{F}{\simeq }$ $^{\Box }\mathcal{C}$ and $%
M\otimes _{^{\ast }\mathcal{C}}$ $^{\Box }\mathcal{C}\overset{\varpi _{M}}{%
\simeq }M$ for every $M\in \mathcal{M}^{\mathcal{C}}.$
\end{proposition}

\begin{Beweis}
Consider for every $M\in \mathcal{M}^{\mathcal{C}}$ the commutative diagram 
\begin{equation*}
\xymatrix{M \otimes_{\mathcal{^*{C}}} Q \otimes_B A \ar[rr]^(.55){\xi _{M}
\otimes id_A} \ar[d]_{id_M \otimes F} & & M^{co{\mathcal{C}}} \otimes_B A
\ar[d]^{\Psi_M} \\ M \otimes_{^*{\mathcal{C}}} {}^{\Box }{\mathcal{C}}
\ar[rr]_{\varpi_M} & & M }
\end{equation*}

Assume $\func{Im}(F)=$ $^{\Box }\mathcal{C}$ and $\varpi _{M}$ to be
surjective for every $M\in \mathcal{M}^{\mathcal{C}}.$ Then $\Psi _{M}$ is
obviously surjective. Let $K=\mathrm{Ke}(\Psi _{M}).$ Since $\Psi _{M}$ is a
morphism in $\mathcal{M}^{\mathcal{C}}\simeq \sigma \lbrack \mathcal{C}%
_{^{\ast }\mathcal{C}}]$ we have $K\in \mathcal{M}^{\mathcal{C}},$ hence $%
\Psi _{K}:K^{co\mathcal{C}}\otimes _{B}A\longrightarrow K$ is surjective. By
Theorem \ref{surj} we have $K\otimes _{^{\ast }\mathcal{C}}Q\overset{\xi _{K}%
}{\simeq }K^{co\mathcal{C}}$ and $A\otimes _{^{\ast }\mathcal{C}}Q\overset{%
\xi _{A}}{\simeq }B,$ hence 
\begin{equation*}
K^{co\mathcal{C}}\simeq K\otimes _{^{\ast }\mathcal{C}}Q=\mathrm{Ke}(\Psi
_{M})\otimes _{^{\ast }\mathcal{C}}Q=\mathrm{Ke}(\Psi _{M}\otimes _{^{\ast }%
\mathcal{C}}id_{Q})=\mathrm{Ke}(id_{M^{co\mathcal{C}}})=0,
\end{equation*}
i.e. $\Psi _{M}$ is bijective. By corollary \ref{co=x} $\Phi _{N}$ is
bijective for every $N\in \mathcal{M}_{B}.$ Consequently $\mathcal{M}^{%
\mathcal{C}}$ satisfies the strong structure theorem.

On the other hand, assume that $\mathcal{M}^{\mathcal{C}}$ satisfies the
strong structure theorem. Note that $F$ is the adjunction of $\Psi _{^{\Box }%
\mathcal{C}},$ hence $Q\otimes _{B}A\overset{F}{\simeq }$ $^{\Box }\mathcal{C%
}$ and consequently $\varpi _{M}$ is also bijective for every $M\in \mathcal{%
M}^{\mathcal{C}}$ by the commutativity of the above diagram.$\blacksquare $
\end{Beweis}

\begin{remarks}
Assume $_{A}\mathcal{C}$ to be locally projective.

\begin{enumerate}
\item  $\varpi _{A}:A\otimes _{^{\ast }\mathcal{C}}$ $^{\Box }\mathcal{C}%
\longrightarrow A$ is surjective iff there exists $\widehat{g}\in $ $^{\Box }%
\mathcal{C}$ with $\widehat{g}(\mathbf{x})=1_{A}.$ To prove this assume
first that $\varpi _{A}$ is surjective. Then there exist $%
\{(a_{i},g_{i})\}_{i=1}^{k}\subset A\times $ $^{\Box }\mathcal{C},$ such
that $\sum\limits_{i=1}^{k}a_{i}\leftharpoonup g_{i}=1_{A}.$ Set $\widehat{g}%
:=\sum\limits_{i=1}^{k}a_{i}g_{i}\in $ $^{\Box }\mathcal{C}.$ Then $\widehat{%
g}(\mathbf{x})=(\sum\limits_{i=1}^{k}a_{i}g_{i})(\mathbf{x}%
)=\sum\limits_{i=1}^{k}g_{i}(\mathbf{x}a_{i})=\sum\limits_{i=1}^{k}a_{i}%
\leftharpoonup g_{i}=1_{A}.$ On the other hand, assume there exists $%
\widehat{g}\in $ $^{\Box }\mathcal{C}$ with $\widehat{g}(\mathbf{x})=1_{A}.$
Then for every $a\in A$ we have $1_{A}\leftharpoonup (\widehat{g}a)=(%
\widehat{g}a)(\mathbf{x})=\widehat{g}(\mathbf{x})a=a,$ i.e. $\varpi _{A}$ is
surjective.

\item  Assume $\varpi _{A}$ to be surjective. If $\Psi _{M}$ is surjective
for $M\in \mathcal{M}^{\mathcal{C}},$ then $\varpi _{M}$ is surjective,
since 
\begin{equation*}
\varpi _{M}\circ (\Psi _{M}\otimes _{^{\ast }\mathcal{C}}id_{^{\Box }%
\mathcal{C}})=\Psi _{M}\circ (id_{M^{\mathbf{co}\mathcal{C}}}\otimes
_{B}\varpi _{A}).
\end{equation*}
\end{enumerate}
\end{remarks}

\qquad

\begin{theorem}
\label{fin-gen}Assume $_{A}\mathcal{C}$ to be f.g. and projective. Then the
following are equivalent:

\begin{enumerate}
\item  $\mathcal{M}^{\mathcal{C}}$ satisfies the weak structure theorem;

\item  $_{B}A$ is flat and $A/B$ is $\mathcal{C}$-Galois;

\item  $_{B}A$ is flat and $\beta ^{\prime }:=\Psi _{\mathcal{C}}^{\prime }$
is an isomorphism;

\item  $_{B}A$ is flat and for every $A$-generated $M\in \mathcal{M}^{%
\mathcal{C}}=\mathcal{M}_{^{\ast }\mathcal{C}},$ the $\mathcal{C}$-colinear
morphism $\Psi _{M}^{\prime }$ is bijective;

\item  for every $M\in \mathcal{M}^{\mathcal{C}},$ the $\mathcal{C}$%
-colinear morphism $\Psi _{M}^{\prime }$ is bijective;

\item  $_{B}A$ is flat, $\mathcal{M}_{^{\ast }\mathcal{C}}=\sigma \lbrack
A_{^{\ast }\mathcal{C}}]$ and $\mathrm{Hom}_{-^{\ast }\mathcal{C}}(A,-):%
\mathrm{Gen}(A_{^{\ast }\mathcal{C}})\longrightarrow \mathcal{M}_{B}$ is
full faithful;

\item  $\mathrm{Hom}_{-^{\ast }\mathcal{C}}(A,-):\mathcal{M}_{^{\ast }%
\mathcal{C}}\longrightarrow \mathcal{M}_{B}$ is faithful;

\item  $A_{^{\ast }\mathcal{C}}$ is a generator;

\item  $F:Q\otimes _{B}A\longrightarrow $ $^{\ast }\mathcal{C}$ is
surjective \emph{(}bijective\emph{)};

\item  (a)\ $Q_{B}$ is f.g. and projective;

(b) $\Omega :A\longrightarrow \mathrm{Hom}_{-B}(Q,B),$ $a\mapsto \lbrack
q\mapsto a\leftharpoonup q]$ is a bimodule isomorphism;

(c) $_{^{\ast }\mathcal{C}}Q$ is faithful;

\item  (a)$\;_{B}A$ is f.g. and projective;

(b)\ $\Lambda :$ $^{\ast }\mathcal{C}\longrightarrow \mathrm{End}%
(_{B}A)^{op},$ $g\mapsto \lbrack a\mapsto a\leftharpoonup g]$ is a ring
isomorphism.
\end{enumerate}
\end{theorem}

\begin{Beweis}
The result follows by Theorems \ref{C -finite}, \ref{gen} and the fact that
in case $_{A}\mathcal{C}$ is f.g. and projective $\mathcal{M}^{\mathcal{C}}=%
\mathcal{M}_{^{\ast }\mathcal{C}}=\sigma \lbrack \mathcal{C}_{^{\ast }%
\mathcal{C}}].\blacksquare $
\end{Beweis}

\begin{theorem}
\label{morita}\emph{(Morita, e.g. \cite[4.1.3, 4.3]{Fai81}, \cite[2.6]{MZ97}%
). }Let $\mathcal{R}$ be a ring, $P$ a right $\mathcal{R}$-module, $\mathcal{%
S}:=\mathrm{End}(P_{\mathcal{R}})$ and $P^{\ast }:=\mathrm{Hom}_{\mathcal{R}%
}(P,\mathcal{R}).$

\begin{enumerate}
\item  The following are equivalent:

\begin{enumerate}
\item  $P_{\mathcal{R}}$ is a generator;

\item  $_{\mathcal{S}}P$ is f.g. projective and $\mathcal{R}\simeq \mathrm{%
End}(_{\mathcal{S}}P)^{op}$ canonically.
\end{enumerate}

\item  The following are equivalent:

\begin{enumerate}
\item  $P_{\mathcal{R}}$ is a faithful quasiprogenerator and $_{\mathcal{S}%
}P $ is finitely generated;

\item  $P_{\mathcal{R}}$ is a progenerator;

\item  $_{\mathcal{S}}P$ is a progenerator and $P_{\mathcal{R}}$ is
faithfully balanced;

\item  $P_{\mathcal{R}}$ and $_{\mathcal{S}}P$ are generators;

\item  $P_{\mathcal{R}}$ and $_{\mathcal{S}}P$ are f.g. and projective;

\item  $\mathrm{Hom}_{-\mathcal{R}}(P,-):\mathcal{M}_{\mathcal{R}%
}\longrightarrow \mathcal{M}_{\mathcal{S}}$ is an equivalence with inverse $%
\mathrm{Hom}_{-\mathcal{S}}(P^{\ast },-);$

\item  $-\otimes _{\mathcal{R}}P^{\ast }:\mathcal{M}_{\mathcal{R}%
}\longrightarrow \mathcal{M}_{\mathcal{S}}$ is an equivalence with inverse $%
-\otimes _{\mathcal{S}}P.$
\end{enumerate}
\end{enumerate}
\end{theorem}

As a consequence of Theorems \ref{prog} and \ref{morita} we get

\begin{theorem}
\label{fin-prog}Assume $_{A}\mathcal{C}$ to be f.g. and projective. Then the
following are equivalent:

\begin{enumerate}
\item  $\mathcal{M}^{\mathcal{C}}$ satisfies the strong structure theorem;

\item  $_{B}A$ is faithfully flat and $A/B$ is $\mathcal{C}$-Galois;

\item  $_{B}A$ is faithfully flat and $\beta ^{\prime }:=\Psi _{\mathcal{C}%
}^{\prime }$ is bijective;

\item  $_{B}A$ is faithfully flat and for every $M\in \sigma \lbrack
A_{^{\ast }\mathcal{C}}],$ the map $\Psi _{M}^{\prime }$ is bijective;

\item  $A_{^{\ast }\mathcal{C}}$ is quasiprojective and generates each of
its submodules, $_{B}A$ is a weak generator and $\mathcal{M}_{^{\ast }%
\mathcal{C}}=\sigma \lbrack A_{^{\ast }\mathcal{C}}];$

\item  $A_{^{\ast }\mathcal{C}}$ is a quasiprogenerator and $\mathcal{M}%
_{^{\ast }\mathcal{C}}=\sigma \lbrack A_{^{\ast }\mathcal{C}}];$

\item  $_{B}A$ is a weak generator, $\Psi _{M}^{\prime }$ is an isomorphism
for every $M\in \mathrm{Gen}(A_{^{\ast }\mathcal{C}})$ and $\mathcal{M}%
_{^{\ast }\mathcal{C}}=\sigma \lbrack A_{^{\ast }\mathcal{C}}];$

\item  $A_{^{\ast }\mathcal{C}}$ is a faithful quasiprogenerator and $_{B}A$
is finitely generated;

\item  $_{B}A$ is a progenerator and $A_{^{\ast }\mathcal{C}}$ is faithfully
balanced;

\item  $\mathrm{Hom}_{-^{\ast }\mathcal{C}}(A,-):\mathcal{M}_{^{\ast }%
\mathcal{C}}\longrightarrow \mathcal{M}_{B}$ is an equivalence with inverse $%
\mathrm{Hom}_{-B}(Q,-);$

\item  $-\otimes _{^{\ast }\mathcal{C}}Q:\mathcal{M}_{^{\ast }\mathcal{C}%
}\longrightarrow \mathcal{M}_{B}$ is an equivalence with inverse $-\otimes
_{B}A;$

\item  $A_{^{\ast }\mathcal{C}}$ and $_{B}A$ are generators;

\item  $A_{^{\ast }\mathcal{C}}$ and $_{B}A$ are f.g. and projective;

\item  $A_{^{\ast }\mathcal{C}}$ is a progenerator.
\end{enumerate}
\end{theorem}

\section{Cleft $C$-Galois Extensions}

In what follows $R$ is a \emph{commutative }ring with $1_{R}\neq 0_{R}$ and $%
\mathcal{M}_{R}$ is the category of $R$-(bi)modules. For an $R$-coalgebra $%
(C,\Delta _{C},\varepsilon _{C})$ and an $R$-algebra $(A,\mu _{A},\eta _{A})$
we consider $(\mathrm{Hom}_{R}(C,A),\star ):=\mathrm{Hom}_{R}(C,A)$ as an $R$%
-algebra with the \emph{usual convolution product }$(f\star g)(c):=\sum
f(c_{1})g(c_{2})$ and unity $\eta _{A}\circ \varepsilon _{C}.$ The unadorned 
$-\otimes -$ means $-\otimes _{R}-.$

\begin{punto}
\textbf{Entwined modules.}\label{ent-str}\ A right-right \textbf{entwining
structure }$(A,C,\psi )$ over $R$ consists of an $R$-algebra $(A,\mu
_{A},\eta _{A}),$ an $R$-coalgebra $(C,\Delta _{C},\varepsilon _{C})$ and an 
$R$-linear map 
\begin{equation*}
\psi :C\otimes _{R}A\longrightarrow A\otimes _{R}C,\text{ }c\otimes a\mapsto
\sum a_{\psi }\otimes c^{\psi },
\end{equation*}
such that 
\begin{equation*}
\begin{tabular}{llllll}
$\sum (a\widetilde{a})_{\psi }\otimes c^{\psi }$ & $=$ & $\sum a_{\psi }%
\widetilde{a}_{\Psi }\otimes c^{\psi \Psi },$ & $\sum (1_{A})_{\psi }\otimes
c^{\psi }$ & $=$ & $1_{A}\otimes c,$ \\ 
$\sum a_{\psi }\otimes \Delta _{C}(c^{\psi })$ & $=$ & $\sum a_{\psi \Psi
}\otimes c_{1}^{\Psi }\otimes c_{2}^{\psi },$ & $\sum a_{\psi }\varepsilon
_{C}(c^{\psi })$ & $=$ & $\varepsilon _{C}(c)a.$%
\end{tabular}
\end{equation*}
\end{punto}

\begin{punto}
Let $(A,C,\psi )$ be a right-right entwining structure. An \textbf{entwined
module }corresponding to $(A,C,\psi )$ is a right $A$-module, which is also
a right $C$-comodule through $\varrho _{M},$ such that 
\begin{equation*}
\varrho _{M}(ma)=\sum m_{<0>}a_{\psi }\otimes m_{<1>}^{\psi }\text{ for all }%
m\in M\text{ and }a\in A.
\end{equation*}
The category of right-right entwined modules and $A$-linear $C$-colinear
morphisms is denoted by $\mathcal{M}_{A}^{C}(\psi ).$ For $M,N\in \mathcal{M}%
_{A}^{C}(\psi )$ we denote by $\mathrm{Hom}_{A}^{C}(M,N)$ the set of $A$%
-linear $C$-colinear morphisms from $M$ to $N.$ With $\#_{\psi }^{op}(C,A):=%
\mathrm{Hom}_{R}(C,A),$ we denote the $A$-ring with $(af)(c)=\sum a_{\psi
}f(c^{\psi }),$ $(fa)(c)=f(c)a,$ multiplication $(f\cdot g)(c)=\sum
f(c_{2})_{\psi }g(c_{1}^{\psi })$ and unity $\eta _{A}\circ \varepsilon _{C}$
(see \cite[Lemma 3.3]{Abu03}).

Entwined modules were introduced by T. Brzezi\'{n}ski and S. Majid \cite
{BM98} as a generalization of the Doi-Koppinen modules presented in \cite
{Doi92} and \cite{Kop95}. By a remark of M. Takeuchi (e.g. \cite[Proposition
2.2]{Brz02}), we have an $A$-coring structure on $\mathcal{C}:=A\otimes
_{R}C,$ where $\mathcal{C}$ is an $A$-bimodule through $a(\widetilde{a}%
\otimes c):=a\widetilde{a}\otimes c,$ $(\widetilde{a}\otimes c)a:=\sum 
\widetilde{a}a_{\psi }\otimes c^{\psi }$ and has comultiplication 
\begin{equation*}
\Delta _{\mathcal{C}}:A\otimes _{R}C\longrightarrow (A\otimes _{R}C)\otimes
_{A}(A\otimes _{R}C),\text{ }a\otimes c\mapsto \sum (a\otimes c_{1})\otimes
_{A}(1_{A}\otimes c_{2})
\end{equation*}
and counity $\varepsilon _{\mathcal{C}}:=id_{A}\otimes \varepsilon _{C}.$
Moreover $\mathcal{M}_{A}^{C}(\psi )\simeq \mathcal{M}^{\mathcal{C}},$ $%
\#_{\psi }^{op}(C,A)\simeq $ $^{\ast }\mathcal{C}$ as $A$-rings and $_{A}%
\mathcal{C}$ is flat (resp. f.g., projective), if $_{R}C$ is so (e.g. \cite
{Abu03}).
\end{punto}

\qquad Inspired by \cite[3.1]{Doi94} we make the following definition:

\begin{punto}
\label{s-rat}Let $(A,C,\psi )$ be a right-right entwining structure over $R$
and consider the corresponding $A$-coring $\mathcal{C}:=A\otimes _{R}C.$ We
say that $(A,C,\psi )$ satisfies the \textbf{left }$\alpha $\textbf{%
-condition}, if for every right $A$-module $M$ the following map is
injective 
\begin{equation*}
\alpha _{M}^{\psi }:M\otimes _{R}C\longrightarrow \mathrm{Hom}_{R}(\#_{\psi
}^{op}(C,A),M),\text{ }m\otimes c\mapsto \lbrack f\mapsto mf(c)]
\end{equation*}
(equivalently, if $_{A}\mathcal{C}$ is locally projective).

Let $M$ be a right $\#_{\psi }^{op}(C,A)$-module $M$ and consider the
canonical map $\rho _{M}:M\longrightarrow \mathrm{Hom}_{R}(\#_{\psi
}^{op}(C,A),M).$ Set $\mathrm{Rat}^{C}(M_{\#_{\psi }^{op}(C,A)}):=(\rho
_{M}^{\psi })^{-1}(M\otimes _{R}C).$ We call $M$ $\#$\textbf{-rational}, if $%
\mathrm{Rat}^{C}(M_{\#_{\psi }^{op}(C,A)})=M$ and set $\varrho _{M}:=(\alpha
_{M}^{\psi })^{-1}\circ \rho _{M}.$ The category of $\#$-rational right $%
\#_{\psi }^{op}(C,A)$-modules will be denoted by $\mathrm{Rat}^{C}(\mathcal{M%
}_{\#_{\psi }^{op}(C,A)}).$
\end{punto}

\begin{theorem}
\label{ent-sg}\emph{(\cite[Theorem 3.10]{Abu03})}\ Let $(A,C,\psi )$ be a
right-right entwining structure and consider the corresponding $A$-coring $%
\mathcal{C}:=A\otimes _{R}C.$

\begin{enumerate}
\item  If $_{R}C$ is flat, then $\mathcal{M}_{A}^{C}(\psi )$ is a
Grothendieck category with enough injective objects.

\item  If $_{R}C$ is locally projective \emph{(}resp. f.g. and projective%
\emph{)}, then 
\begin{equation}
\mathcal{M}_{A}^{C}(\psi )\simeq \mathrm{Rat}^{C}(\mathcal{M}_{\#_{\psi
}^{op}(C,A)})\simeq \sigma \lbrack (A\otimes _{R}C)_{\#_{\psi }^{op}(C,A)}]\;%
\text{\emph{(}resp. }\mathcal{M}_{A}^{C}(\psi )\simeq \mathcal{M}_{\#_{\psi
}^{op}(C,A)}\text{\emph{)}.}  \label{iso-sm}
\end{equation}
\end{enumerate}
\end{theorem}

In what follows we fix a right-right entwining structure $(A,C,\psi )$ with $%
\mathcal{C}:=A\otimes _{R}C$ the corresponding $A$-coring and assume that $%
A\in \mathcal{M}_{A}^{C}(\psi )\simeq \mathcal{M}^{\mathcal{C}}$ with 
\begin{equation*}
\varrho _{A}:A\longrightarrow A\otimes _{R}C,\text{ }a\mapsto \sum
a_{<0>}\otimes a_{<1>}=\sum 1_{<0>}a_{\psi }\otimes 1_{<1>}^{\psi }.
\end{equation*}
Then $\sum 1_{<0>}\otimes 1_{<1>}\in \mathcal{C}$ is a group-like element
and 
\begin{equation*}
Q\simeq \{q\in \mathrm{Hom}_{R}(C,A)\mid \sum q(c_{2})_{\psi }\otimes
c_{1}^{\psi }=\sum q(c)1_{<0>}\otimes 1_{<1>}\text{ for all }c\in C\}.
\end{equation*}
For every $M\in \mathcal{M}_{A}^{C}(\psi ),$ we set 
\begin{equation*}
M^{co\mathcal{C}}:=\{m\in M\mid \sum m_{<0>}\otimes m_{<1>}=\sum
m1_{<0>}\otimes 1_{<1>}\}.
\end{equation*}
Moreover we set $B:=A^{co\mathcal{C}}.$

\begin{remark}
Let $x\in C$ be a group-like element. For every right $C$-comodule $M$ we
put $M^{coC}:=\{m\in M\mid \varrho _{M}(m)=m\otimes x\}.$ If $\varrho
_{A}(1_{A})=1_{A}\otimes x,$ then we have $M^{coC}=M^{co\mathcal{C}}$ for
every $M\in \mathcal{M}_{A}^{C}(\psi ).$
\end{remark}

\qquad By \cite[Corollaries 3.4, 3.7]{Brz99} $-\otimes _{R}^{c}A:\mathcal{M}%
^{C}\longrightarrow \mathcal{M}_{A}^{C}(\psi )$ is a functor, which is left
adjoint to the forgetful functor. Here, for every $N\in \mathcal{M}^{C},$ we
consider the canonical right $A$-module $N\otimes _{R}^{c}A:=N\otimes _{R}A$
with the $C$-coaction $n\otimes a\mapsto \sum n_{<0>}\otimes a_{\psi
}\otimes n_{<1>}^{\psi }.$

\begin{proposition}
\label{sp}Let $R$ be a QF ring and assume $C$ be right semiperfect. Let $%
_{R}C$ to be locally projective (projective) and put $C^{\Box }:=\mathrm{Rat}%
(_{^{\ast }C}C^{\ast }).$

\begin{enumerate}
\item  The following are equivalent:

\begin{enumerate}
\item  $A$ is a generator in $\mathcal{M}_{A}^{C}(\psi );$

\item  $A$ generates $C^{\Box }\otimes _{R}^{c}A$ in $\mathcal{M}%
_{A}^{C}(\psi );$

\item  the map $\Psi _{C^{\Box }\otimes _{R}^{c}A}^{\prime }:\mathrm{Hom}%
_{A}^{C}(A,C^{\Box }\otimes _{R}^{c}A)\otimes _{B}A\longrightarrow C^{\Box
}\otimes _{R}^{c}A$ is surjective \emph{(}bijective\emph{)}.
\end{enumerate}

\item  The following are equivalent:

\begin{enumerate}
\item  $A$ is a progenerator in $\mathcal{M}_{A}^{C}(\psi );$

\item  $\Psi _{C^{\Box }\otimes _{R}^{c}A}^{\prime }$ is surjective \emph{(}%
bijective\emph{)} and $_{B}A$ is a weak generator.
\end{enumerate}
\end{enumerate}
\end{proposition}

\begin{Beweis}
By \cite[2.6]{MTW01} $C^{\Box }$ is a generator in $\mathcal{M}^{C},$ hence $%
C^{\Box }\otimes _{R}^{c}A$ is a generator in $\mathcal{M}_{A}^{C}(\psi )$
by the functorial isomorphism $\mathrm{Hom}_{A}^{C}(C^{\Box }\otimes
_{R}^{c}A,M)\simeq \mathrm{Hom}^{C}(C^{\Box },M)$ for every $M\in \mathcal{M}%
_{A}^{C}(\psi ).$

\begin{enumerate}
\item  The assertions follow form the note above and Theorem \ref{gen}.

\item  (a) $\Rightarrow $ (b) follows by Theorem \ref{prog}.

(b)\ $\Rightarrow $ (a).\ By the note above $C^{\Box }\otimes _{R}^{c}A$ is
a generator in $\mathcal{M}_{A}^{C}(\psi ),$ and the surjectivity of $\Psi
_{C^{\Box }\otimes _{R}^{c}A}^{\prime }$ makes $A$ a generator in $\mathcal{M%
}_{A}^{C}(\psi ).$ So $_{B}A$ is flat by Theorem \ref{gen}. The weak
generator property makes $_{B}A$ faithfully flat and we are done by Theorem 
\ref{prog}.$\blacksquare $
\end{enumerate}
\end{Beweis}

\begin{definition}
A \textbf{(total) integral} for $C$ is a $C$-colinear morphism $\lambda
:C\longrightarrow A$ (with $\sum 1_{<0>}\lambda (1_{<1>})=1_{A}$). We call
the ring extension $A/B$ \textbf{cleft, }if there exists a $\star $\textbf{%
-invertible integral}. We say $A$ has the \textbf{right normal basis property%
}, if there exists a left $B$-linear right $C$-colinear isomorphism $A\simeq
B\otimes _{R}C.$
\end{definition}

\begin{lemma}
\label{co-Q}Let $\lambda \in \mathrm{Hom}_{R}(C,A)$ be $\star $-invertible
with inverse $\overline{\lambda }.$ Then:

\begin{enumerate}
\item  $\lambda \in \mathrm{Hom}^{C}(C,A)$ iff $\overline{\lambda }\in Q.$

\item  If $\varrho (a)=\sum a_{\psi }\otimes x^{\psi }$ for some group-like
element $x\in C,$ then there exists $\widehat{\lambda }\in Q,$ such that $%
\sum 1_{<0>}\widehat{\lambda }(1_{<1>})=\widehat{\lambda }(x)=1_{A}$ \emph{(}%
in this case $C$ admits a total integral, namely the $\star $-inverse of $%
\widehat{\lambda }$\emph{)}.
\end{enumerate}
\end{lemma}

\begin{Beweis}
Let $\lambda \in \mathrm{Hom}_{R}(C,A)$ be $\star $-invertible with inverse $%
\overline{\lambda }.$

\begin{enumerate}
\item  If $\overline{\lambda }\in Q,$ then we have for all $c\in C:$%
\begin{equation*}
\begin{tabular}{lllll}
$\sum \lambda (c_{1})\otimes c_{2}$ & $=$ & $\sum \lambda (c_{1})1_{\psi
}\otimes c_{2}^{\psi }$ &  &  \\ 
& $=$ & $\sum \lambda (c_{1})\varepsilon (c_{3})1_{\psi }\otimes c_{2}^{\psi
}$ &  &  \\ 
& $=$ & $\sum \lambda (c_{1})(\overline{\lambda }(c_{3})\lambda
(c_{4}))_{\psi }\otimes c_{2}^{\psi }$ &  &  \\ 
& $=$ & $\sum \lambda (c_{1})\overline{\lambda }(c_{3})_{\psi }\lambda
(c_{4})_{\Psi }\otimes c_{2}^{\psi \Psi }$ &  &  \\ 
& $=$ & $\sum \lambda (c_{1})\overline{\lambda }(c_{22})_{\psi }\lambda
(c_{3})_{\Psi }\otimes c_{21}^{\psi \Psi }$ &  &  \\ 
& $=$ & $\sum \lambda (c_{1})\overline{\lambda }(c_{2})1_{<0>}\lambda
(c_{3})_{\Psi }\otimes 1_{<1>}^{\Psi }$ &  &  \\ 
& $=$ & $\sum 1_{<0>}\lambda (c)_{\Psi }\otimes 1_{<1>}^{\Psi }$ &  &  \\ 
& $=$ & $\sum \lambda (c)_{<0>}\otimes \lambda (c)_{<1>}$ &  & 
\end{tabular}
\end{equation*}
i.e. $\lambda \in \mathrm{Hom}^{C}(C,A).$ On the other hand, if $\lambda \in 
\mathrm{Hom}^{C}(C,A),$ then we have for all $c\in C:$%
\begin{equation*}
\begin{tabular}{lllll}
$\sum \overline{\lambda }(c_{2})_{\psi }\otimes c_{1}^{\psi }$ & $=$ & $\sum 
\overline{\lambda }(c_{1})\lambda (c_{2})\overline{\lambda }(c_{4})_{\psi
}\otimes c_{3}^{\psi }$ &  &  \\ 
& $=$ & $\sum \overline{\lambda }(c_{1})\lambda (c_{2})_{<0>}\overline{%
\lambda }(c_{3})_{\psi }\otimes \lambda (c_{2})_{<1>}^{\psi }$ &  &  \\ 
& $=$ & $\sum \overline{\lambda }(c_{1})1_{<0>}\lambda (c_{2})_{\psi }%
\overline{\lambda }(c_{3})_{\Psi }\otimes 1_{<1>}^{\psi \Psi }$ &  &  \\ 
& $=$ & $\sum \overline{\lambda }(c_{1})1_{<0>}(\lambda (c_{2})\overline{%
\lambda }(c_{3}))_{\psi }\otimes 1_{<1>}^{\psi }$ &  &  \\ 
& $=$ & $\sum \overline{\lambda }(c)1_{<0>}1_{\psi }\otimes 1_{<1>}^{\psi }$
&  &  \\ 
& $=$ & $\sum \overline{\lambda }(c)1_{<0>}\otimes 1_{<1>},$ &  & 
\end{tabular}
\end{equation*}
i.e. $\overline{\lambda }\in Q.$

\item  Assume $\varrho (a)=\sum a_{\psi }\otimes x^{\psi }$ for some
group-like element $x\in C.$ Let $\lambda \in \mathrm{Hom}^{C}(C,A)$ with $%
\overline{\lambda }\in Q$ (see (1)). Then $\widehat{\lambda }:=\overline{%
\lambda }\lambda (x)\in Q,$ since $\lambda (x)\in B,$ and moreover $\sum
1_{<0>}\widehat{\lambda }(1_{<1>})=\widehat{\lambda }(x)=\overline{\lambda }%
(x)\lambda (x)=(\overline{\lambda }\star \lambda )(x)=\varepsilon
_{C}(x)1_{A}=1_{A}.\blacksquare $
\end{enumerate}
\end{Beweis}

\begin{proposition}
\label{cleft}Assume $A/B$ to be cleft.

\begin{enumerate}
\item  $\mathcal{M}_{A}^{C}(\psi )$ satisfies the weak structure theorem $%
\emph{(}$in particular $A/B$ is $\mathcal{C}$-Galois\emph{)}.

\item  For every $M\in \mathcal{M}_{A}^{C}(\psi ),$ the $C$-colinear
morphism 
\begin{equation*}
\gamma _{M}:M\longrightarrow M^{co\mathcal{C}}\otimes _{R}C,\text{ }m\mapsto
\sum m_{<0>}\overline{\lambda }\otimes m_{<1>}
\end{equation*}
is an isomorphism.

\item  $A$ has the right normal basis property.

\item  If $_{R}C$ is faithfully flat, then $\mathcal{M}_{A}^{C}(\psi )$
satisfies the strong structure theorem.
\end{enumerate}
\end{proposition}

\begin{Beweis}
Assume there exists a $\star $-invertible $\lambda \in \mathrm{Hom}^{C}(C,A)$
with inverse $\overline{\lambda }\in Q$ (see Lemma \ref{co-Q} (1)).

\begin{enumerate}
\item  Let $M\in \mathcal{M}_{A}^{C}(\psi )$ and consider 
\begin{equation*}
\widetilde{\Psi }_{M}:M\longrightarrow M^{co\mathcal{C}}\otimes _{B}A,\text{ 
}m\mapsto \sum m_{<0>}\overline{\lambda }\otimes \lambda (m_{<1>}).
\end{equation*}
Then we have for all $n\in M^{co\mathcal{C}},$ $m\in M$ and $a\in A:$%
\begin{equation*}
\begin{tabular}{lllll}
$(\widetilde{\Psi }_{M}\circ \Psi _{M})(n\otimes a)$ & $=$ & $\widetilde{%
\Psi }_{M}(na)$ &  &  \\ 
& $=$ & $\sum (na_{<0>})\overline{\lambda }\otimes _{B}\lambda (a_{<1>})$ & 
&  \\ 
& $=$ & $\sum na_{<0><0>}\overline{\lambda }(a_{<0><1>})\otimes _{B}\lambda
(a_{<1>})$ &  &  \\ 
& $=$ & $\sum n\otimes _{B}a_{<0><0>}\overline{\lambda }(a_{<0><1>})\lambda
(a_{<1>})$ &  &  \\ 
& $=$ & $\sum n\otimes _{B}a_{<0>}\overline{\lambda }(a_{<1>1})\lambda
(a_{<1>2})$ &  &  \\ 
& $=$ & $n\otimes _{B}a$ &  & 
\end{tabular}
\end{equation*}
and 
\begin{equation*}
\begin{tabular}{lllll}
$(\Psi _{M}\circ \widetilde{\Psi }_{M})(m)$ & $=$ & $\sum (m_{<0>}\overline{%
\lambda })\lambda (m_{<1>})$ &  &  \\ 
& $=$ & $\sum m_{<0><0>}\overline{\lambda }(m_{<0><1>})\lambda (m_{<1>})$ & 
&  \\ 
& $=$ & $\sum m_{<0>}\overline{\lambda }(m_{<1>1})\lambda (m_{<1>2})$ &  & 
\\ 
& $=$ & $\sum m_{<0>}\varepsilon _{C}(m_{<1>})1_{A}$ &  &  \\ 
& $=$ & $m.$ &  & 
\end{tabular}
\end{equation*}

\item  For every $M\in \mathcal{M}_{A}^{C}(\psi ),$ $\gamma _{M}$ is
bijective with inverse 
\begin{equation*}
\widetilde{\gamma }_{M}:M^{co\mathcal{C}}\otimes _{R}C\longrightarrow M,%
\text{ }n\otimes c\mapsto n\lambda (c).
\end{equation*}
In fact we have for all $m\in M,$ $n\in M^{co\mathcal{C}}$ and $c\in C:$%
\begin{equation*}
\begin{tabular}{lllll}
$(\widetilde{\gamma }_{M}\circ \gamma _{M})(m)$ & $=$ & $\sum (m_{<0>}%
\overline{\lambda })\lambda (m_{<1>})$ &  &  \\ 
& $=$ & $\sum m_{<0><0>}\overline{\lambda }(m_{<0><1>})\lambda (m_{<1>})$ & 
&  \\ 
& $=$ & $\sum m_{<0>}\overline{\lambda }(m_{<1>1})\lambda (m_{<1>2})$ &  & 
\\ 
& $=$ & $\sum m_{<0>}\varepsilon _{C}(m_{<1>})$ &  &  \\ 
& $=$ & $m$ &  & 
\end{tabular}
\end{equation*}
and 
\begin{equation*}
\begin{tabular}{lllll}
$(\gamma _{M}\circ \widetilde{\gamma }_{M})(n\otimes _{B}c)$ & $=$ & $\sum
(n\lambda (c))_{<0>}\overline{\lambda }\otimes (n\lambda (c))_{<0>}$ &  & 
\\ 
& $=$ & $\sum (n\lambda (c)_{<0>})\overline{\lambda }\otimes \lambda
(c)_{<1>}$ &  &  \\ 
& $=$ & $\sum (n\lambda (c_{1}))\overline{\lambda }\otimes c_{2}$ &  &  \\ 
& $=$ & $\sum n\lambda (c_{1})_{<0>}\overline{\lambda }(\lambda
(c_{1})_{<1>})\otimes c_{2}$ &  &  \\ 
& $=$ & $\sum n\lambda (c_{11})\overline{\lambda }(c_{12})\otimes c_{2}$ & 
&  \\ 
& $=$ & $n\otimes c.$ &  & 
\end{tabular}
\end{equation*}

\item  By (2)\ the left $B$-linear right $C$-colinear map 
\begin{equation*}
\gamma _{A}:A\longrightarrow B\otimes _{R}C,a\mapsto \sum
a_{<0>}\leftharpoonup \overline{\lambda }\otimes a_{<1>}
\end{equation*}
is an isomorphism with inverse $b\otimes c\mapsto b\lambda (c).$

\item  Assume $_{R}C$ to be faithfully flat. By (3) $A\simeq B\otimes _{R}C$
as left $B$-modules, hence $_{B}A$ is faithfully flat. By (1) $A/B$ is $%
\mathcal{C}$-Glaois and we are done by Theorem \ref{prog}.$\blacksquare $
\end{enumerate}
\end{Beweis}

\begin{theorem}
\label{main}The following statements are equivalent:

\begin{enumerate}
\item  $A/B$ is cleft;

\item  $\mathcal{M}_{A}^{C}(\psi )$ satisfies the weak structure theorem and 
$A$ has the right normal basis property;

\item  $A/B$ is $\mathcal{C}$-Galois and $A$ has the right normal basis
property;

\item  $\Lambda :\#_{\psi }^{op}(C,A)\simeq \mathrm{End}(_{B}A)^{op},$ $%
g\mapsto \lbrack a\mapsto a\leftharpoonup $ $g]$ is a ring isomorphism and $%
A $ has the right normal basis property.

If moreover $_{R}C$ is faithfully flat, then \emph{(1)-(4)} are equivalent to

\item  $\mathcal{M}_{A}^{C}(\psi )$ satisfies the strong structure theorem
and $A$ has the right normal basis property.
\end{enumerate}
\end{theorem}

\begin{Beweis}
(1)\ $\Rightarrow $ (2). This follows by Proposition \ref{cleft}.

(2) $\Rightarrow $ (3). By assumption $\beta :=\Psi _{A\otimes _{R}C}$ is an
isomorphism.

(3) $\Rightarrow $ (4). By assumption $A\otimes _{B}A\simeq A\otimes _{R}C$
as left $A$-modules, hence we have the canonical isomorphisms 
\begin{equation*}
\begin{tabular}{lllll}
$\#_{\psi }^{op}(C,A)$ & $\simeq $ & $\mathrm{Hom}_{A-}(A\otimes _{R}C,A)$ & 
$\simeq $ & $\mathrm{Hom}_{A-}(A\otimes _{B}A,A)$ \\ 
& $\simeq $ & $\mathrm{Hom}_{B-}(A,\mathrm{End}(_{A}A))$ & $\simeq $ & $%
\mathrm{End}(_{B}A).$%
\end{tabular}
\end{equation*}

(4)\ $\Rightarrow $ (1). Assume $\theta :B\otimes _{R}C\longrightarrow A$ to
be a left $B$-linear right $C$-colinear isomorphism and consider the right $%
C $-colinear morphism $\lambda :C\longrightarrow A,$ $c\mapsto \theta
(1_{A}\otimes c)$ and the left $B$-linear morphism $\delta :=(id\otimes
\varepsilon _{C})\circ \theta ^{-1}:A\longrightarrow B.$ Define $\overline{%
\lambda }:=\Lambda ^{-1}(\delta )\in \#_{\psi }^{op}(C,A).$ Then we have for
all $c\in C:$%
\begin{equation*}
\begin{tabular}{lllll}
$\sum \lambda (c_{1})\overline{\lambda }(c_{2})$ & $=$ & $\sum \lambda
(c)_{<0>}\overline{\lambda }(\lambda (c)_{<1>})$ & $=$ & $\lambda
(c)\leftharpoonup \overline{\lambda }$ \\ 
& $=$ & $\delta (\lambda (c))$ & $=$ & $((id\otimes \varepsilon _{C})\circ
\theta ^{-1})(\lambda (c))$ \\ 
& $=$ & $((id\otimes \varepsilon _{C})\circ \theta ^{-1})(\theta
(1_{A}\otimes c))$ & $=$ & $\varepsilon _{C}(c)1_{A}.$%
\end{tabular}
\end{equation*}
On the other hand we have for all $a\in A:$%
\begin{equation*}
\begin{tabular}{lllll}
$\Lambda (\overline{\lambda }\star \lambda )(a)$ & $=$ & $a\leftharpoonup (%
\overline{\lambda }\star \lambda )$ & $=$ & $\sum a_{<0>}(\overline{\lambda }%
\star \lambda )(a_{<1>})$ \\ 
& $=$ & $\sum a_{<0>}\overline{\lambda }(a_{<1>1})\lambda (a_{<1>2})$ & $=$
& $\sum a_{<0><0>}\overline{\lambda }(a_{<0><1>})\lambda (a_{<1>})$ \\ 
& $=$ & $\sum (a_{<0>}\leftharpoonup \overline{\lambda })\lambda (a_{<1>})$
& $=$ & $\sum (a_{<0>}\leftharpoonup \Lambda ^{-1}(\delta ))\lambda
(a_{<1>}) $ \\ 
& $=$ & $\sum \delta (a_{<0>})\lambda (a_{<1>})$ & $=$ & $\sum \delta
(a_{<0>})\theta (1_{A}\otimes a_{<1>})$ \\ 
& $=$ & $\sum \theta (\delta (a_{<0>})\otimes a_{<1>})$ & $=$ & $\theta
(\theta ^{-1}(a))=a,$%
\end{tabular}
\end{equation*}
hence $\overline{\lambda }\star \lambda =\eta _{A}\circ \varepsilon _{C}.$

Now assume $_{R}C$ to be faithfully flat. Then (1) $\Rightarrow $ (5)
follows by Proposition \ref{cleft} (4) and we are done.$\blacksquare $
\end{Beweis}

\qquad The following result deals with the special case $\varrho (a)=\sum
a_{\psi }\otimes x^{\psi },$ for some group-like element $x\in C.$ In this
case we obtain the equivalent statements (1)-(5) in Theorem \ref{main}
without any assumptions on $C.$

\begin{theorem}
\label{x-case}Assume that $\varrho (a)=\sum a_{\psi }\otimes x^{\psi }$ for
some group-like element $x\in C.$ The following statements are equivalent:

\begin{enumerate}
\item  $A/B$ is cleft;

\item  $\mathcal{M}_{A}^{C}(\psi )$ satisfies the strong structure theorem
and $A$ has the right normal basis property;

\item  $\mathcal{M}_{A}^{C}(\psi )$ satisfies the weak structure theorem and 
$A$ has the right normal basis property;

\item  $A/B$ is $\mathcal{C}$-Galois and $A$ has the right normal basis
property;

\item  $\Lambda :\#_{\psi }^{op}(C,A)\simeq \mathrm{End}(_{B}A)^{op},$ $%
g\mapsto \lbrack a\mapsto a\leftharpoonup $ $g]$ is a ring isomorphism and $%
A $ has the right normal basis property.
\end{enumerate}
\end{theorem}

\begin{Beweis}
By Theorem \ref{main} it remains to prove that $\Phi _{N}$ is an isomorphism
for every $N\in \mathcal{M}_{B},$ if $A/B$ is cleft. But in our special case
there exists by Lemma \ref{co-Q} some $\widehat{\lambda }\in Q$ with $\sum
1_{<0>}\widehat{\lambda }(1_{<1>})=1_{A}$ and we are done by Corollary \ref
{co=x} (2).$\blacksquare $
\end{Beweis}

\begin{remark}
Let $(H,A,C)$ be a right-right resp. a left-right\ Doi-Koppinen structure.
Then $(A,C,\psi )$ is a right-right entwining structure with 
\begin{equation*}
\psi :C\otimes _{R}A\longrightarrow A\otimes _{R}C,\text{ }c\otimes a\mapsto
\sum a_{<0>}\otimes ca_{<1>}
\end{equation*}
resp. a left-right entwining structure with 
\begin{equation*}
\psi :A\otimes _{R}C\longrightarrow A\otimes _{R}C,\text{ }a\otimes c\mapsto
\sum a_{<0>}\otimes a_{<1>}c.
\end{equation*}
If $x$ is a group-like element of $C,$ then $A\in \mathcal{M}(H)_{A}^{C}$
with $\varrho (a):=\sum a_{<0>}\otimes xa_{<1>}$ (resp. $\varrho (a)=\sum
a_{<0>}\otimes a_{<1>}x$) and we get \cite[Theorem 1.5]{DM92} (resp. 
\cite[Theorem 2.5]{Doi94}) as special cases of Theorem \ref{x-case}.
\end{remark}

\end{document}